\documentclass[11pt]{article}

\usepackage[margin=1in]{geometry}
\usepackage{amsmath} 
\usepackage{amssymb} 
\usepackage{amsbsy} 
\usepackage{amsthm} 
\usepackage{amsfonts}
\usepackage{thmtools}
\usepackage{thm-restate}
\usepackage{bbm}
\usepackage{url}
\usepackage{booktabs}
\usepackage{mathrsfs}
\usepackage{multirow}
\usepackage{chngpage}
\usepackage{subcaption}
\usepackage{float}
\usepackage{subfiles}
\usepackage[math]{cellspace}
\usepackage[font={small}]{caption}
\usepackage{tikz, pgfplots}
\usepackage{algpseudocode}
\usepackage{graphicx}
\usepackage{algorithm}
\usepackage{fullpage}
\usepackage{chngcntr}
\usepackage[square,sort,comma,numbers]{natbib}
\usepackage{soul}
\usepackage{mathtools}
\usepackage{array}
\usepackage{arydshln}
\usepackage[letterspace=-85]{microtype}
\usepackage{nicefrac}
\usepackage{hyperref}
\usepackage{xcolor}

\usetikzlibrary{shapes.geometric,arrows,arrows.meta,positioning,calc}
\usetikzlibrary{patterns}
\usetikzlibrary{decorations.pathmorphing,backgrounds,fit,matrix,decorations.pathreplacing}

\newcommand*{\pfstart}{\begin{proof}}
\newcommand*{\pfend}{\end{proof}}
\newcommand{\blue}[1]{\textcolor{black}{#1}}

\newcommand{\one}{\mathbf{1}}

\tikzstyle{vertex} = [circle, minimum size=0.1cm, inner sep=0pt, draw=black, fill=black]
\tikzstyle{circ} = [circle, minimum width=0.5mm, inner sep=0pt,draw,fill]
\tikzstyle{hcirc} = [circle,minimum width=5mm, inner sep=0pt,draw]
\tikzstyle{bcirc} = [circle, minimum width=1.5mm, inner sep=0pt,draw,fill]
\tikzstyle{bhcirc} = [circle, minimum width=1.5mm, inner sep=0pt,draw, dotted ]
\tikzstyle{ept} = [circle,minimum width=0mm, inner sep=0pt, white]
\tikzstyle{txt} = [text width=1.3cm,draw,rounded corners=3pt]
\tikzstyle{ncirc} = [circle,draw=black, inner sep=1pt]

\DeclarePairedDelimiter\iprod{\langle}{\rangle}

\algrenewcommand\algorithmicindent{1em}

\newcommand{\SectionLabel}[1]{\Statex\hspace*{-\algorithmicindent}\textbf{#1}}

\theoremstyle{remark}

\usepackage{bm}
\theoremstyle{remark}
\newtheorem*{claim*}{Claim}
\theoremstyle{remark}
\newtheorem*{remark*}{Remark}
\theoremstyle{remark}
\newtheorem{remark}{Remark}

\theoremstyle{plain}

\theoremstyle{plain}

\theoremstyle{definition}

\theoremstyle{definition}
\newtheorem*{assumption*}{Assumption}

\theoremstyle{definition}
\newtheorem*{example*}{Example}

\renewcommand\thmcontinues[1]{Cont.}

\theoremstyle{definition}

\theoremstyle{plain}

\theoremstyle{plain}
\newtheorem{corollary}{Corollary}

\newcommand{\R}{\mathbb{R}}

\newcommand{\hL}{\hat{\mathcal L}}
\DeclareMathOperator{\argmin}{argmin}

\newcommand{\nn}{\nonumber}

\newcommand{\dom}{\mathrm{dom}\,}

\hypersetup{
    colorlinks=true,
    citecolor=black,
    linkcolor=black,
    filecolor=black,
    urlcolor=black,
}

\pgfplotsset{compat=1.15}
\tikzset{snake it/.style={decorate, decoration={snake, amplitude=.3mm, segment length=2mm}}}

\allowdisplaybreaks

\newcommand{\ABSTRACT}[1]{\begin{abstract}#1\end{abstract}}

\newcommand{\KEYWORDS}[1]{\par\medskip\noindent\textbf{Keywords: }#1\par}

\begin{document}

\title{GPU-Accelerated First-Order Method with Randomized Sampling for Binary Integer Programs}
\author{
    Ningji Wei
    \thanks{
    Department of Industrial, Manufacturing and Systems Engineering, Texas Tech University, \texttt{ningji.wei@ttu.edu}. This work was partially supported by AFOSR grant FA9550-26-1-B345 and NSF grant 2548574.}
    \and
    Jiaming Liang
    \thanks{
    Goergen Institute for Data Science and Artificial Intelligence (GIDS-AI) and Department of Computer Science, University of Rochester, \texttt{jiaming.liang@rochester.edu}. This work was partially supported by GIDS-AI seed funding and AFOSR grant FA9550-25-1-0182.}
  }
\date{}
\maketitle

\ABSTRACT{%
We present a scalable, GPU-accelerated algorithmic framework for large binary integer programs that operates end-to-end with minimal synchronization overhead. The proposed method combines a first-order routine that guides search in the continuous relaxation with a randomized, feasibility-aware sampling module that generates batched binary candidates. 
% We establish near-stationary-point guarantees for the first-order component, together with probabilistic bounds on the feasibility and solution quality of sampled candidates, without providing global optimality certificates.
\blue{We establish a residual convergence guarantee for the first-order component in the convex setting and a finite-horizon primal-dual residual estimate in the general nonconvex setting, together with probabilistic bounds on the feasibility and solution quality of sampled candidates, without providing global optimality certificates.}
A key modeling consideration under this framework concerns the treatment of equality constraints to enhance sampling effectiveness. We address this through total-unimodular reformulations and tailored sampling schemes designed to improve feasibility rates.
We evaluate the framework on standard benchmark classes, including set cover, knapsack, max-cut, three-dimensional assignment, and facility location. On small- to medium-scale instances, state-of-the-art exact solvers remain stronger overall; nevertheless, our method produces high-quality incumbents within short runtimes. On larger instances, the proposed approach achieves substantially shorter runtimes while often delivering solution quality comparable to exact solvers under the same time limit. Overall, the framework complements exact methods by enabling scalable, parallel search in domains where problem size and response time are primary considerations.
}%

% \FUNDING{The first author was partially supported by AFOSR grant FA9550-26-1-B345 and NSF grant 2548574. The second author was partially supported by GIDS-AI seed funding and AFOSR grant FA9550-25-1-0182.}

%Supplemental Material:
%Data Ethics & Reproducibility Note:

% Sample
%\KEYWORDS{Stochastic programming, Decision support,Uncertainty, Disaster response, Optimization}

% Fill in data. If unknown, outcomment the field
%\HISTORY{Received: Month DD, YYYY; Accepted: Month DD, YYYY; Published Online: Month DD, YYYY}

\section{Introduction}
Over the past decade, advances in GPU-accelerated algorithms have enabled substantial speedups for core machine learning workloads \citep{dean2012large}.
These gains have also motivated GPU-based developments in large-scale continuous optimization, including linear programming \citep{lu2023cupdlp} and other first-order methods \citep{o2021operator,stellato2020osqp,shin2024accelerating}. 
It is therefore natural to ask whether similar benefits can be realized for large-scale mixed integer programming (MIP).
Despite progress on specific problem classes \citep{boyer2012solving,suri2012scalable} and GPU-accelerated primal heuristics that assist classic solvers \citep{ccorduk2025gpu}, a central challenge remains: the tree-based search at the core of branch-and-bound is inherently sequential and synchronization-heavy, which limits the effectiveness of fully parallel execution on modern GPU architectures.
In this paper, we take a step toward addressing this gap by focusing on the following class of binary integer programs (BIPs):
\begin{subequations}
  \label{eq:bip}
  \begin{align}
  \min_{x \in \{0,1\}^n} &~\iprod{x, Qx} + \iprod{c, x} \\
                         &~ A x \geq b \label{eq:bip01}\\
                         &~ B x  = d, \label{eq:bip02}
  \end{align}
\end{subequations}
where $A \in \mathbb{R}^{m_1 \times n}$ and $B \in \mathbb{R}^{m_2 \times n}$ encode the inequality and equality constraints, and $Q \in \mathbb{R}^{n \times n}$ is assumed, without loss of generality, to be symmetric but not necessarily positive semidefinite (PSD).

This problem class encompasses many canonical combinatorial optimization problems (e.g., set cover, knapsack) and is NP-hard in general. The area has been extensively studied: a large body of literature spans problem-specific algorithms \citep{rendl2010solving,applegate2011traveling}, exact formulations and approximations \citep{pardalos1994maximum,caprara2000approximation}, and general-purpose solvers \citep{bixby2012brief}. These developments have expanded the reach of exact solvers to tackle instance sizes previously impractical \citep{koch2022progress}, enabling decision support across domains such as defense, transportation, and communications.

Meanwhile, the scale and time requirements of many modern applications continue to increase, with growing needs to obtain high-quality feasible solutions quickly on large-scale instances, further motivating scalable approaches for large BIPs under tight decision windows.
We propose the \emph{GPU-accelerated first-order method with randomized sampling} (GFORS), an end-to-end framework designed to operate efficiently on massively parallel architectures.
We establish \blue{a finite-horizon primal-dual residual estimate} for the underlying first-order routine, together with probabilistic guarantees on feasibility and solution quality for the sampling scheme, albeit without providing global optimality certificates.
The GFORS framework consists of three components:

\begin{itemize}
  \item \emph{First-Order Method for Region Exploration:} 
    We apply a first-order method to \eqref{eq:bip} that steers fractional solutions toward promising regions. \blue{Section~\ref{sec:foalgo} establishes for the implemented \emph{primal-dual hybrid gradient} (PDHG) variant \citep{chambolle2011first} a residual convergence result in the convex setting and a finite-horizon residual estimate in the nonconvex setting.}

\item \emph{Randomized Sampling for Candidate Generation:} 
Each iteration of the first-order method yields a fractional solution $x \in [0,1]^n$, which admits a natural probabilistic interpretation (Section~\ref{sec:rsalgo}). This enables a sampling subroutine guided by the iterate trajectory of the first-order component. In Section~\ref{sec:sanalysis}, we will derive probability bounds on feasibility and objective quality for this subroutine.

\item \emph{GPU Acceleration for Scalability:} 
Both components are designed for GPU execution, enabling the framework to run with minimal CPU-GPU synchronization and to fully leverage the scalability of modern GPU architectures.
\end{itemize}

To strengthen the framework and improve sampling efficiency, we introduce additional techniques such as total-unimodular reformulation and customized sampling design. We evaluate the resulting framework on a range of classic BIPs, including set cover, knapsack, max-cut, three-dimensional assignment, and facility location. An aggregated overview is presented in Figure~\ref{fig:perf-both} (see Section~\ref{sec:agg}): on small- to medium-scale instances, Gurobi \blue{ and cuOpt} achieve stronger performance, yet our method still delivers high-quality solutions within runtimes of seconds; for large-scale instances, however, our method achieves substantially shorter runtimes in general, 
while its solution quality is often comparable to the exact solvers under the prescribed time limit.
This is \emph{not} a head-to-head algorithmic comparison; rather, the results illustrate how GFORS can complement exact methods by offering scalable performance on large BIP instances with short decision windows.

\subsection{Related Literature}

\paragraph{Exact solvers.} 
Each state-of-the-art exact MIP solver (e.g., Gurobi, CPLEX, SCIP, Xpress, and CBC) deploys a portfolio of techniques, including strong presolve \citep{achterberg2020presolve,gamrath2015progress}, cutting-plane separation \citep{ceria1998cutting}, local-search heuristics \citep{fischetti2003local,danna2005exploring,fischetti2005feasibility}, and parallel branch-and-bound/branch-and-cut \citep{applegate2011traveling}, to deliver certified optimality on large, heterogeneous instances. 
These components are combined within a parallelized tree search with adaptive cut/heuristic scheduling, representing the prevailing standard for exact MILP/MIQP/BIP.

% For comprehensive overviews, we refer the reader to \citet{bolusani2024scip,gurobi2025docs,bixby1999mip}.

\paragraph{First-order methods.}
First-order methods such as gradient descent, proximal algorithms \citep{moreau1965proximite,rockafellar1976monotone}, and PDHG \citep{chambolle2011first} are widely used in large-scale convex optimization due to their low per-iteration cost and amenability to parallelization. They have been extensively applied in continuous domains such as signal processing \citep{combettes2011proximal}, imaging \citep{chambolle2011first}, and, more recently, large-scale linear programming \citep{lu2023cupdlp}. Their use in MIP, however, remains nascent, since the loss of convexity arising from integrality limits direct applicability.

\paragraph{Randomized Rounding.}
Randomized rounding \citep{raghavan1987randomized} was originally introduced as a post-processing technique: after solving the linear programming relaxation of a binary linear program, the fractional solution is rounded randomly to produce binary solutions. The approach has since been adapted to many problem classes, including the traveling salesman problem \citep{gharan2011randomized}, the largest simplex problem \citep{nikolov2015randomized}, virtual network embedding \citep{rost2019virtual}, and machine scheduling \citep{schulz2002scheduling}. More recently, dependent rounding variants introduce designed rounding correlations, still starting from a fractional optimum, to better exploit problem structure \citep{bertsimas1999dependent}. We refer the reader to \citet{srinivasan1999approximation} for a comprehensive survey. 
Although closely related, our randomized sampling differs in three key respects: (i) it is embedded within a first-order iterative method (rather than applied after solving a linear program); (ii) it draws samples repeatedly from a dynamic distribution that evolves with the algorithm's updates across iterations; and (iii) it is GPU-oriented, enabling massively parallel, batched sampling and evaluation each iteration.

\paragraph{GPU-accelerated algorithms \blue{for continuous optimization}.} 
GPUs are now central to large-scale machine learning, where massive data and models benefit from highly parallel, vectorized first-order methods with low per-iteration cost, enabling training deep neural networks via gradient-based algorithms \citep{dean2012large}. Beyond machine learning, GPU-accelerated solvers exist for linear, conic, and nonlinear optimization, including cuPDLP \citep{lu2023cupdlp}, SCS \citep{o2021operator}, OSQP \citep{stellato2020osqp}, and MadNLP.jl \citep{shin2024accelerating}. Notably, cuPDLP demonstrates substantial performance gains on large-scale linear programs. 

\paragraph{\blue{GPU-accelerated algorithms for discrete optimization.}}
\blue{
Compared with continuous optimization, relatively few studies have explored GPU-based algorithms for MIPs. Early work largely focused on specific problem classes \citep{boyer2012solving,suri2012scalable}, with a broader survey provided by \citet{boyer2017gpu}. More recently, GPU-accelerated first-order methods have been integrated into discrete optimization in several ways. For problems with strong continuous relaxations, scalable gradient-based methods have been developed for \(K\)-means clustering and \(k\)-sparse generalized linear models \citep{de2026power,liu2025scalable,liu2026gpu}. GPU-accelerated primal heuristics have also been used to support CPU-based branch-and-bound (B\&B), including fix-and-propagate guided by inexact PDLP solutions \citep{kempke2026fix} and massively parallel repair, local-search, and feasibility heuristics in cuOpt and CHAP \citep{ccorduk2025gpu,tjusila2026chap}. More recently, batched first-order methods have been developed for LP subproblems arising within exact MIP procedures \citep{blin2026batched}. These approaches mainly use GPU computation to augment the conventional CPU-based B\&B framework. In contrast, GFORS is designed as a GPU-native method, where a first-order algorithm evolves a probability distribution over binary candidates and iterative sampling is used to update feasible incumbents.
}

\paragraph{Quadratic unconstrained binary optimization (QUBO).}
QUBO addresses problems of the form $\min_{x\in\{0,1\}^n} \iprod{x, Qx} + \iprod{c, x}$ with \emph{no constraints}, which are equivalent to Ising formulations and thus attractive for specialized hardware and sampling-based methods \citep{kochenberger2014unconstrained,lucas2014ising}. Leading approaches include classical metaheuristics, exact branch-and-bound variants, and quantum/quantum-inspired algorithms, e.g., quantum annealing \citep{kadowaki1998quantum}, that naturally return batches of samples per evaluation. Despite a shared probabilistic perspective, the mechanisms differ: quantum-based methods sample from fixed, parameterized states over unconstrained encodings, whereas GFORS performs adaptive, GPU-batched sampling guided by first-order information and uses dual updates to explicitly drive constraint feasibility.

% Throughout the paper, for a matrix $A$ and index sets $J$ (rows) and $I$ (columns), we write $A_{J}$ for the submatrix of rows indexed by $J$, $A_{\cdot I}$ for the submatrix of columns indexed by $I$, and $A_{JI}$ for the submatrix restricted to rows $J$ and columns $I$. For a single index $j \in J$ and $i \in I$, $A_{j}$, $A_{\cdot i}$, and $A_{ji}$ denote the associated row, column, and entry. For any function $f$ differentiable at $y$, we use $\ell_f(\cdot;y)$ as its affine approximation at $y$.

\subsection{Contributions}
We summarize our main contributions as follows.
\begin{itemize}

\item \emph{A GPU-native first-order \& sampling framework for large-scale BIPs.}
We formulate the BIP solution as the iterative refinement of a product distribution over binary solutions, where the continuous relaxation is solved by PDHG and integrated with massively parallel, batched sampling in an end-to-end GPU execution framework.

  \item \emph{Theoretical guarantees for the relaxation and sampling mechanisms.}
  % We establish near-stationarity convergence guarantees for PDHG on solving the weakly convex relaxation and derive probabilistic bounds on the feasibility and objective quality of sampled solutions for general BIPs.
    \blue{We establish an ${\cal O}(N^{-1/2})$ best-iterate residual convergence guarantee for PDHG in the convex setting and a finite-horizon primal-dual residual estimate in the general nonconvex setting, together with probabilistic bounds on the feasibility and objective quality of sampled solutions for~\eqref{eq:bip}.}

  \item \emph{Sampling enhancement under equality constraints.}
  We incorporate techniques such as total-unimodular reformulations and customized sampling schemes into the framework to restore sampling effectiveness in equality-constrained settings.

  \item \emph{Large-scale computational validation.}
  Across multiple classical BIP classes, we demonstrate strong time-to-incumbent performance relative to the exact solvers Gurobi \blue{ and cuOpt} on large-scale instances.
\end{itemize}

Throughout the paper, for a matrix $A$ and index sets $J$ (rows) and $I$ (columns), we denote by $A_J$ the submatrix consisting of rows in $J$, by $A_{\cdot I}$ the submatrix consisting of columns in $I$, and by $A_{JI}$ the submatrix restricted to both rows $J$ and columns $I$. For individual indices $j \in J$ and $i \in I$, $A_j$, $A_{\cdot i}$, and $A_{ji}$ refer to the corresponding row, column, and entry, respectively. For any function $f$ differentiable at $y$, the affine function $\ell_f(\cdot; y)$ denotes its first-order approximation at $y$. We use $(x_*, y_*)$ to denote a saddle point, and $x^*$ to denote an optimal solution to an optimization problem. The remainder of the paper is organized as follows. Section~\ref{sec:framework} presents the main design of the GFORS framework. Section~\ref{sec:implement} describes additional implementation details. Section~\ref{sec:expr} reports both high-level and problem-specific computational results. Finally, Section~\ref{sec:conc} concludes with further discussions. To streamline the exposition, we defer all the proofs to Appendix~\ref{sec:proofs}.

% The code repository is available at \url{https://github.com/tidues/GFORS}.

\section{GFORS Framework}
\label{sec:framework}
This section introduces the main design of GFORS for solving BIPs of the form \eqref{eq:bip}. To enable GPU acceleration, we reformulate this problem into the following saddle-point form with an additional fractional penalty term. To ease the notation, let $y=(u,v)$ be the concatenated dual variables for \eqref{eq:bip01} and \eqref{eq:bip02}, $\one$ denotes the all-ones vector, and $K:=-[A; B]$ and $r := (b; d)$ denote the row-wise concatenation of the respective matrices and vectors.
\begin{equation}
  \label{eq:saddle}  
\begin{aligned}
  \min_{x} \max_{y} &~ \left\{\Phi(x; u, v):=\hat{\mathcal L}(x; y) + h_1(x) - h_2(y)\right\}\\
  \text{where} &~ \hat{\mathcal L}(x; y):= \iprod{x, Qx} + \iprod{c, x} + \iprod{y, Kx + r} + \rho \iprod{x, \one-x},\\ 
               &~ h_1(x):=I_{[0,1]^n}(x),\\
               &~ h_2(y):=I_{\mathbb R_+^{m_1}\times \mathbb R^{m_2}}(y).
\end{aligned}
\end{equation}
The function $\hat{\mathcal L}$ denotes the Lagrangian of the linear relaxation of \eqref{eq:bip} along with the fractional penalty term $\rho \langle x, \one-x\rangle$, and $I_{\mathcal X}(x)$ is the indicator function that equals $0$ when $x \in \mathcal X$ and $+\infty$ otherwise.
\blue{Note that $-h_2(y)$ enforces the constrained maximization over the dual feasible set in the extended-real-valued formulation.}
For fixed $\rho$ and $y$, the function $\hat {\mathcal L}(\cdot;y)$ is weakly convex and $L$-smooth with $L:=2(\|Q\| + \rho)$.

\begin{figure}[t]
\centering

%==================== (a) Landscape over [0,1]^2 ====================
\begin{subfigure}[t]{0.49\linewidth}
\centering
% Original concept01.tex design with a slightly larger usable plot area.
\begin{tikzpicture}[trim axis left, trim axis right]
  \pgfmathdeclarefunction{Lhat}{2}{%
    \pgfmathparse{((#1-0.78)^2 + (#2-0.18)^2)
                 + 0.08*sin(deg(8*#1))*cos(deg(7*#2))}%
  }
  \def\zlift{1}

  \begin{axis}[
    width=1.10\linewidth,
    height=0.825\linewidth,
    view={80}{20},
    xmin=0, xmax=1,
    ymin=0, ymax=1,
    zmin=0,
    xlabel={},
    ylabel={},
    zlabel={$\hat{\mathcal L}(\cdot;y)$},
    ticks=none,
    axis lines=left,
    axis line style={line width=0.6pt, -{Stealth[length=2mm,width=1.3mm]}},
    grid=none,
    clip=false,
    zlabel style={yshift=0pt},
  ]

  \addplot3[
    mesh,
    samples=28, samples y=28,
    domain=0:1, y domain=0:1,
    very thin,
    draw=black,
  ]
  {Lhat(x,y) + \zlift};

  \addplot3[black, line width=0.6pt]
    coordinates {(0,0,0) (1,0,0) (1,1,0) (0,1,0) (0,0,0)};

  \addplot3[only marks, mark=o, mark size=1.6pt, black]
    coordinates {(0,0,0) (0,1,0) (1,0,0) (1,1,0)};
  \node[anchor=north east, xshift=2pt, yshift=0pt]
    at (axis cs:0,0,0) {$(0,0)$};
  \node[anchor=south east, xshift=4pt, yshift=0pt]
    at (axis cs:0,1,0) {$(0,1)$};
  \node[anchor=north west, xshift=-2pt, yshift=-2pt]
    at (axis cs:0.8,-0.15,0) {$(1,0)$};
  \node[anchor=south west, xshift=2pt, yshift=0pt]
    at (axis cs:1.45,0.8,0) {$(1,1)$};

  \addplot3[black, thick, -stealth]
    coordinates {
      (0.50,0.50,0.02)
      (0.56,0.62,0.02)
      (0.68,0.58,0.02)
      (0.62,0.44,0.02)
      (0.75,0.46,0.02)
      (0.70,0.30,0.02)
      (0.83,0.33,0.02)
      (0.78,0.20,0.02)
      (0.90,0.18,0.02)
      (0.86,0.12,0.02)
    };

  \addplot3[only marks, mark=*, mark size=1.2pt, black]
    coordinates {(0.68,0.58,0.02) (0.78,0.20,0.02)};
  \node[anchor=south west] at (axis cs:1.05,0.5,0) {$x_{k_1}$};
  \node[anchor=south west] at (axis cs:0.78,0.1,0) {$x_{k_2}$};
  \end{axis}
\end{tikzpicture}
\caption{$\hat{\mathcal L}(\cdot;y)$ over $[0,1]^2$ with the solution trajectory.}
\label{fig:landscape_2d}
\end{subfigure}
\hfill
%==================== (b) Distributions over {00,01,10,11} ====================
\begin{subfigure}[t]{0.49\linewidth}
\centering
% Original concept02.tex design with tighter margins and a larger plot area.
\begin{tikzpicture}[trim axis left, trim axis right]
  \begin{axis}[
    width=1.10\linewidth,
    height=0.77\linewidth,
    ybar,
    bar width=8pt,
    ymin=0, ymax=0.75,
    ylabel={Probability},
    ylabel style={yshift=-27pt},
    symbolic x coords={00,01,10,11},
    xtick=data,
    xticklabels={{$(0,0)$},{$(0,1)$},{$(1,0)$},{$(1,1)$}},
    enlarge x limits=0.20,
    axis lines*=left,
    ytick=\empty,
    clip=true,
    legend style={
      draw=none,
      at={(0.32,0.96)},
      anchor=north east,
      legend cell align=left,
    },
    axis line style={line width=0.6pt, -{Stealth[length=2mm,width=1.3mm]}},
  ]

  \addplot+[
    draw=black,
    fill=white,
    pattern=north east lines,
    pattern color=black,
    bar shift=-5pt,
  ] coordinates {
    (00,0.1344) (01,0.1856) (10,0.2856) (11,0.3944)
  };

  \addplot+[
    draw=black,
    fill=white,
    bar shift=10pt,
  ] coordinates {
    (00,0.1760) (01,0.0440) (10,0.6240) (11,0.1560)
  };

  \legend{$\pi(x_{k_1})$,$\pi(x_{k_2})$}
  \end{axis}
\end{tikzpicture}
\caption{Distributions over $\{0,1\}^2$ induced by $x_{k_1}$ and  $x_{k_2}$.}
\label{fig:pi_bars}
\end{subfigure}

\caption{Illustration of the objective $\hat{\mathcal L}(\cdot, y):[0,1]^2\to\mathbb{R}$ guiding the first-order iterates, the resulting solution trajectory in the relaxed space $[0,1]^2$, and the sampling distributions induced by two selected iterates along this trajectory.}
\label{fig:two_panel_landscape_pi}
\end{figure}

Building on this formulation, GFORS couples two modules: a first-order method that updates a fractional iterate in the relaxation, and a randomized sampler that draws nearby binary candidates. The first-order dynamics steer $x$ toward regions that achieve low objective values, while sampling translates the current fractional state into a batch of discrete solutions for evaluation and feasibility checks.
Formally, given $\hat{\mathcal L}(\cdot; y): [0,1]^n \to \mathbb{R}$, we define
\[
\pi:\ [0,1]^n \to \mathcal{P}(\{0,1\}^n),
\]
where $\mathcal{P}(\{0,1\}^n)$ denotes the set of probability distributions on $\{0,1\}^n$. Each $x\in[0,1]^n$ thus parameterizes a distribution $\pi(x)$, from which we sample a binary candidate $z\sim \pi(x)$. Figure~\ref{fig:two_panel_landscape_pi} provides an illustration.

\begin{algorithm}[t]
\caption{GFORS Framework}\label{algo:gfors}
\small
\begin{algorithmic}
  \State Parameters: binary penalty $\rho$; sampling interval $k_{\text{int}}$; sampling rounds $k_r$; batch size $k_b$
\SectionLabel{CPU:}
\State \textsc{TUReformulate}()   \Comment{perform TU reformulation if necessary (Section~\ref{sec:tureform})}
\State \textsc{Preprocess}()    \Comment{normalize model parameters (Section~\ref{sec:implement})}

\SectionLabel{GPU:}
\State initialization: $x_0 \in [0,1]^n$; $k \gets 0$; $z_{\text{best}} \gets$ \blue{$+\infty$}; $x_{\text{best}} \gets \varnothing$; \text{halt} $\gets \text{false}$
\Comment{\blue{assuming minimization}}
\While{not \text{halt}}
  \State $k \gets k+1$
  \State $\rho \gets \textsc{UpdatePenalty}()$  \Comment{increase penalty by schedule (Section~\ref{sec:implement})}
  \State $x_k \gets \textsc{FirstOrderStep}(x_{k-1}, \rho)$ \Comment{update $x$ by a first-order method (Section~\ref{sec:foalgo})}
  \If{$k \equiv 0 \pmod{k_{\text{int}}}$} \Comment{sampling trigger for efficiency}

  \For{$i \in [k_r]$}
    \State $\bar X^k \gets \textsc{RandSampleStep}(x_k, k_b)$  \Comment{randomized/customized sampling  (Section~\ref{sec:rsalgo})}
    \State $(z_{\text{best}}, x_{\text{best}}) \gets \textsc{EvalBest}(\bar X^k, z_{\text{best}}, x_{\text{best}})$;
  \EndFor
  \EndIf
  \State \text{halt} $\gets \textsc{CheckHalt}()$; \Comment{stopping criteria (Section~\ref{sec:implement})}

\EndWhile

\State $(z_{\text{best}}, x_{\text{best}}) \gets \textsc{EvalBest}(\mathrm{round}(x_k), z_{\text{best}}, x_{\text{best}})$;
\SectionLabel{CPU:}
\State \Return $z_{\text{best}},\ x_{\text{best}}$
\end{algorithmic}
\end{algorithm}

In this setting, the first-order component operates entirely in the relaxed space, using $\hat{\mathcal L}$ to drive progress. For weakly convex objectives, many first-order schemes admit convergence guarantees to stationary points \citep{ghadimi2016accelerated,liang2021average,kong2019complexity,kong2021accelerated,liang2021fista}. We implement a PDHG variant \citep{zhu2008efficient,chambolle2011first}, motivated by its strong empirical performance on GPU-accelerated linear programs \citep{lu2023cupdlp}.
After each first-order update, the current iterate $x_k$ induces a distribution $\pi(x_k)$, which we use for batched sampling and verification. When $x_k$ has high entropy (i.e., remains far from $\{0,1\}^n$), the resulting samples are broadly exploratory; as $x_k$ becomes more polarized, $\pi(x_k)$ increasingly concentrates mass on a smaller set of high-probability binary candidates, enabling targeted refinement. Additionally, we update the binary penalty parameter $\rho$ over iterations to  \blue{heuristically} promote integrality. These mechanisms produce a smooth transition from global exploration to focused local search as the algorithm progresses.

% \blue{which is a practical heuristic (This adaptive schedule of $\rho$ is a practical heuristic and is not covered by the  theory in Section \ref{sec:foalgo}.) [NW: need to refine a bit]}

Since both subroutines are designed to run efficiently on GPU hardware with minimal CPU-GPU synchronization, the entire framework can be executed end-to-end on the GPU once the problem data has been \blue{preprocessed and} transferred. Algorithm~\ref{algo:gfors} summarizes the overall procedure: it begins with CPU-side reformulation and preprocessing, after which the computation proceeds through a GPU-accelerated iterative loop consisting of first-order updates interleaved with batched sampling and verification.

This section therefore focuses on the core functions \textsc{FirstOrderStep} and \textsc{RandSampleStep}, along with the supporting \textsc{TUReformulate} subroutine that enhances sampling when equality constraints are present. Additional subroutines are mainly implementation details and will be deferred to Section~\ref{sec:implement}.

\subsection{First-Order Update Analysis}
\label{sec:foalgo}

The \textsc{FirstOrderStep} subroutine can be instantiated with any first-order method tailored for the saddle-point problem \eqref{eq:saddle}. We adopt the update mechanism of the PDHG algorithm \citep{zhu2008efficient,chambolle2011first}, which has demonstrated strong performance for linear programs and efficient compatibility with GPU acceleration \citep{lu2023cupdlp}.

Algorithm~\ref{algo:foalgo} presents the detailed pseudocode. 
This subsection analyzes its convergence under convex settings 
(e.g., $Q$ is PSD and $\rho = 0$) and estimates the primal-dual stationarity residuals in the more general weakly convex cases. 
\blue{
Throughout the analysis, we treat $\rho$ as a fixed constant. 
Since the $\rho$-scheduling scheme has an upper bound $\rho_{\max}$ (see Section~\ref{sec:implement}), the analysis applies after $\rho$ reaches $\rho_{\max}$, provided that the stated step-size conditions are satisfied.
}
The following definitions will be used to derive these results.

\begin{algorithm}[!t]
  \caption{PDHG Update}\label{algo:foalgo}
\small
\begin{algorithmic}
  \State Parameters: preprocessed input $(Q,K,c,r)$; step size $\tau_1, \tau_2 > 0$
  \Function{FirstOrderStep}{$x_{k-1}, \rho$}
  \State $y_k \gets \Pi_{\mathbb R_+^{m_1}\times \mathbb R^{m_2}}(y_{k-1} + \tau_2 (K \bar x_{k-1}+r))$
  \State $\delta \gets c+\rho \one + K^\top y_{k} + 2Q x_{k-1} - 2\rho x_{k-1}$
  \State $x_k \gets \Pi_{[0,1]^n}(x_{k-1} - \tau_1 \delta)$
  \State $\bar x_k \gets 2x_k - x_{k-1}$
  \State \Return $x_k$
\EndFunction
\end{algorithmic}
\end{algorithm}

\begin{subequations}
\begin{align}
  \text{Gradients:} &~ \nabla_x \hL(x,y) = c+\rho \one + K^\top y + 2 Q x - 2 \rho x,\ \nabla_y \hL(x, y) = Kx + r.\label{eq:dy}\\
  \text{Smoothness:} &~ \|\nabla_x \hL(x_1,y) - \nabla_x \hL(x_2,y)\| \le 2(\|Q\|+\rho)\|x_1-x_2\| = L\|x_1-x_2\|.\label{ineq:L}\\
  \text{PDHG updates:} &~ y_k = \underset{y\in \R^m}\argmin\left\{-\hL(\bar x_{k-1},y) + h_2(y) + \frac{1}{2\tau_2}\|y-y_{k-1}\|^2\right\}, \label{eq:pdhg-y}\\
			&~ x_k = \underset{x\in \R^n}\argmin\left\{\ell_{\hL(\cdot,y_k)}(x;x_{k-1}) + h_1(x) + \frac{1}{2\tau_1}\|x-x_{k-1}\|^2\right\}, \label{eq:pdhg-x} \\
			&~ \bar{x}_k=2x_k-x_{k-1}. \label{eq:pdhg-bx}\\
  \text{Auxiliaries:} &~ x_{-1}=x_0=\bar x_0, \ y^0 = y_0, \ y^k=y_k+\tau_2 K\left(\bar{x}_k- x_k\right).\label{def:yk}
\end{align}
\end{subequations}
In particular, \eqref{ineq:L} ensures that $\hat{\mathcal L}(\cdot, y_k)$ is $L$-smooth, and \eqref{eq:pdhg-y}-\eqref{eq:pdhg-bx} follow the update rules in Algorithm~\ref{algo:foalgo}.

\blue{For Proposition~\ref{prop:convex} and Theorem~\ref{thm:convex}, we assume that the relaxed problem is feasible and satisfies Slater's condition, that is, there exists $\bar x \in \text{ri}([0,1]^n)$ such that $A\bar x >b$ and $B \bar x = d$. 
Under this condition and the convexity of $\hat{\mathcal L}(\cdot, y)$, a saddle point $(x^*,y^*)$ exists.} We first establish a per-iteration optimality characterization for the primal-dual updates.
In particular, each iterate satisfies exact or approximate subgradient inclusions, where the
inexactness is quantified by the linearization error of $\hL(\cdot,y_k)$.

\begin{restatable}{proposition}{propconv}
  \label{prop:convex}
Assume $\hL(\cdot,y_k)$ is convex, $\tau_1 \tau_2 \|K\|^2 \le 1/4$, and $\tau_1 \le 1/(4L)$. Define 
    \begin{equation}\label{def:epsk}
        \varepsilon_k:= \hL(x_k,y_k) - \ell_{\hL(\cdot,y_k)}(x_k;x_{k-1}).
    \end{equation}
    Then, the following statements hold:
    \begin{align}
      \text{(a)}&\ 
            r^y_k:=\frac{y^{k-1} - y^k}{\tau_2} \in \partial \left[-\hL(x_k,\cdot) + h_2\right] (y_k), \ r^x_k:=\frac{x_{k-1} - x_k}{\tau_1} \in \partial_{\varepsilon_k} [\hL(\cdot,y_k) +h_1(\cdot)](x_k);
\label{incl:IPP}
\\
      \text{(b)}&\ \varepsilon_k \le \frac{1}{8 \tau_1} \|x_k-x_{k-1}\|^2, \ \frac{1}{\tau_2}\|y^k-y_k\|^2 \le \frac{1}{4 \tau_1} \|x_k-x_{k-1}\|^2;
      \label{ineq:ek}
      \\
      \text{(c)}&\ \frac{1}{4\tau_1} \sum_{k=1}^N \|x_k-x_{k-1}\|^2+\frac{1}{2\tau_2} \sum_{k=1}^N \|y_k-y^{k-1}\|^2 \le \frac{1}{2\tau_1}\|x_*-x_0\|^2 + \frac{1}{2\tau_2}\|y_*-y_0\|^2.\label{ineq:IPP}
    \end{align}
\end{restatable}

This proposition provides (a) primal-dual residual inclusions, (b) bounds on the
inexactness and iterate differences, and (c) a global summability estimate. As a consequence, we obtain the following iteration-complexity guarantee.

\begin{restatable}{theorem}{thmconv}
  \label{thm:convex}
  Assume $\hat {\mathcal L}(\cdot, y_k)$ is convex. Denote
    \[
    \Delta_0 := \frac{\|x_*-x_0\|^2}{2\tau_1} + \frac{\|y_*-y_0\|^2}{2\tau_2},
    \]
    then for every $N \ge 1$, there exists $k_0 \le N$ such that
    \begin{align}
        &r^y_{k_0} \in \partial \left[-\hL(x_{k_0},\cdot) + h_2\right] (y_{k_0}), \quad r^x_{k_0} \in \partial_{\varepsilon_{k_0}} [\hL(\cdot,y_{k_0}) +h_1(\cdot)](x_{k_0}), \label{con:incl} \\
        &\|r^y_{k_0}\| \le \frac{\sqrt{3\Delta_0}}{\sqrt{\tau_2 N}}, \quad \|r^x_{k_0}\| \le \frac{2\sqrt{\Delta_0}}{\sqrt{\tau_1 N}}, \quad \varepsilon_{k_0} \le \frac{\Delta_0}{2N}. \label{con:ineq}
    \end{align}
\end{restatable}

% -----------primal-------------

% consider
% \[
% \min_x \max_y \iprod{Kx,y} + h_1(x) - h_2(y)
% \]
% equivalent to
% \[
% \min_x h_2^*(Kx) + h_1(x)
% \]

% \begin{align*}
% 	y_k &= \underset{y\in \R^m}\argmin\left\{\iprod{-K \bar x_{k-1},y} + h_2(y) + \frac{1}{2\tau_2}\|y-y_{k-1}\|^2\right\} \\
%     &= \underset{y\in \R^m}\argmax\left\{\iprod{K \bar x_{k-1},y} - h_2^{\tau_2}(y) \right\} \\
%     &= \nabla (h_2^{\tau_2})^*(K \bar x_{k-1}),\\
% 	x_k &= \underset{x\in \R^n}\argmin\left\{ \iprod{Kx, y_k} + h_1(x) + \frac{1}{2\tau_1}\|x-x_{k-1}\|^2\right\} \\
%     &= \underset{x\in \R^n}\argmin\left\{ \iprod{x, K^\top \nabla (h_2^{\tau_2})^*(K \bar x_{k-1})} + h_1(x) + \frac{1}{2\tau_1}\|x-x_{k-1}\|^2\right\} \\
%     &= \underset{x\in \R^n}\argmin\left\{ \ell_{(h_2^{\tau_2})^* \circ K}(x;\bar x_{k-1}) + h_1(x) + \frac{1}{2\tau_1}\|x-x_{k-1}\|^2\right\}, \\
% 	\bar{x}_k &= 2x_k-x_{k-1}. 
% \end{align*}

% \[
% h_2^{\tau_2}(y) = h_2(y) + \frac{1}{2\tau_2}\|y-y_{k-1}\|^2
% \]
% each iteration solves a smooth surrogate (with $y_{k-1}$ implicitly fixed)
% \[
% \min_x (h_2^{\tau_2})^*(Kx) + h_1(x)
% \]

% here, $(h_2^{\tau_2})^*$ is $\tau_2$-smooth

% in the end, we want to solve 
% \[
% \min_x h_2^*(Kx) + h_1(x)
% \]

% -----------primal-------------

\blue{The subsequent Proposition~\ref{prop:sum} and Theorem~\ref{thm:weakconvex} extend the analysis to general nonconvex settings (e.g., $Q$ is non-PSD and/or $\rho>0$) by providing a finite-horizon residual estimate. We first derive a bound on the cumulative iterate differences over the first $N$ iterations.}

%expressed in terms of the objective decrease $\Delta\Phi_N$ and additional controllable terms.

\begin{restatable}{proposition}{propsum}
  \label{prop:sum}
  Assume
\begin{equation}\label{ineq:tau}
    \tau_1 \tau_2 \|K\|^2 \le \frac12, \quad \tau_1\le \frac{1}{2L}.
\end{equation}
Then, the following statement holds:
\begin{equation}\label{ineq:sum}
    \frac{1}{\tau_1} \sum_{k=1}^N \|x_k - x_{k-1}\|^2 + \frac{1}{\tau_2}\sum_{k=1}^N \|y_k - y_{k-1}\|^2 \le 4\Delta\Phi_N + \frac8{\tau_1} \sum_{k=1}^N\|x_{k-1}\|^2 + 16\tau_2 \|r\|^2 N,
\end{equation}
where $\Delta \Phi_N := \Phi(x_0,y_0) - \Phi(x_N,y_N)$.
%Not so desirable, we need to assume $r=0$ and $\sum \|x_k\|^2 \le C$, which might be fine, since $x_k \in [0,1]^n$
\end{restatable}

\blue{Based on the estimate in Proposition~\ref{prop:sum}, Theorem~\ref{thm:weakconvex} provides a finite-horizon bound on the primal-dual residuals defined below, rather than an asymptotic convergence guarantee.}
 
\begin{restatable}{theorem}{thmweakconv}
  \label{thm:weakconvex}
    Denote the primal and dual residuals as
    \[
      \begin{aligned}
        s_k^x &= \frac{x_{k-1} - x_k}{\tau_1} + \nabla_x \hL(x_k,y_k) - \nabla_x \hL(x_{k-1},y_k)\\
        s_k^y &= \frac{y_{k-1} -y_k}{\tau_2} - \nabla_y \hL(x_k,y_k) + \nabla_y \hL(\bar x_{k-1},y_k).
      \end{aligned}
    \]
    % Assume $\Phi(x,y)$ is bounded by $\Phi_*$ from below. 
    Then, for every $N\ge 1$, there exists $k_0\le N$ such that
    \begin{align}
        &s_{k_0}^x \in \nabla_x \hL(x_{k_0},y_{k_0}) + \partial h_1(x_{k_0}), \ s_{k_0}^y \in - \nabla_y \hL(x_{k_0},y_{k_0}) + \partial h_2(y_{k_0}), \label{con:incl1}\\
        &\|s_{k_0}^x\| + \|s_{k_0}^y\| \leq 2 \sqrt{\alpha \beta_N}. \label{bound:xy}
    \end{align}
where $\alpha:=\left( 2\|K\| + \frac{3}{2\tau_1}\right)^2 \tau_1 + \frac{1}{\tau_2}$ and $\beta_N:=\frac{\Delta \Phi_N}{N} + \frac2{\tau_1 N} \sum_{k=1}^N\|x_{k-1}\|^2 + 4\tau_2 \|r\|^2$.
\end{restatable}

\subsection{Randomized Sampling}
\label{sec:rsalgo}

Consider the binary space $\{0,1\}^n$. Any fractional solution $p \in [0,1]^n$ naturally induces a fully independent joint distribution over $\{0,1\}^n$, defined for each $\hat x \in \{0,1\}^n$ as
\[
\pi(p, \hat x) := \Pr(\hat x \mid p) \;=\; \prod_{i \in [n]} p_i^{\hat x_i}(1-p_i)^{1-\hat x_i}.
\]
The associated expectation vector $x_p=\sum_{\hat x \in \{0,1\}^n} \pi(p,\hat x)\hat x$
represents the probabilistic mixture of binary solutions governed by $p$. With these definitions, we obtain the following intuitive proposition.

\begin{restatable}{proposition}{smpinter}
  For every $p \in [0,1]^n$, we have $x_p = p$.
\end{restatable}

Although this equivalence is straightforward, it provides the following useful equivalence
$$\min_{x \in [0,1]^n} f(x) \Longleftrightarrow \min_{p \in [0,1]^n} f(x_p),$$
thereby offering two equivalent interpretations of the first-order optimization process: \blue{it can be viewed either as iteratively refining a fractional solution, or as updating the parameter $p$ of a product distribution over $\{0,1\}^n$ used for subsequent sampling.} Furthermore, if $f$ guarantees that an optimal solution must be binary, then the optimal distribution is indeed fully independent and coincides with the optimal solution $x^*$. Particularly, this interpretation enables the generic Bernoulli sampling subroutine presented in Algorithm~\ref{alg:gsample}. 
% \end{remark}

\begin{algorithm}[b]
\caption{Generic Bernoulli Sampling}
\label{alg:gsample}
\begin{algorithmic}
\Function{RandSampleStep}{$p, k$}
  % \State $\bar X \gets \mathbb{1}\{\,U < \mathbf{1}_k p^\top\,\}$ \;\; where $U \sim \mathrm{Unif}(0,1)^{k \times n}$  \Comment{rows i.i.d.; $X_{ij}\sim \mathrm{Ber}(p_j)$}
\For{$l \in [k], i \in [n]$}
      \State $\bar X_{l, i} \gets \mathrm{Bernoulli}(p_i)$
  \EndFor
  \State \Return $\bar X \in \{0,1\}^{k \times n}$
\EndFunction
\end{algorithmic}
\end{algorithm}

\subsection{Sampling Algorithm Analysis}
\label{sec:sanalysis}
This subsection derives probability bounds that establish the optimality and feasibility guarantees of Algorithm~\ref{alg:gsample}, offering insight into the behavior of GFORS.

\begin{restatable}{theorem}{optprob}
  \label{thm:optprob}
  Let $p \in [0,1]^n$ and $k \in \mathbb{N}$ in Algorithm~\ref{alg:gsample}. 
  Let $f$ and $x^* \in \{0,1\}^n$ be the objective function and an optimizer of~\eqref{eq:bip}, and fix an optimality tolerance $\delta > 0$. 
  Set $\mu := \mathbb{E}_{x \sim p}[\|x - x^*\|_1]$ as the expected 1-norm distance from $x^*$, and let $L_f$ be the Lipschitz constant of $f$ under 1-norm. 
  \blue{Define
  \begin{equation}\label{eq:q-obj}
      q_{\mathrm{obj}} = \exp\!\left( \tfrac{\delta}{L_f}\bigl(1 - \log\tfrac{\delta}{\mu L_f}\bigr) - \mu \right).
  \end{equation}}Then, given $\mu \leq \delta / L_f$, the following bound holds,
  \[
    \Pr\!\left(\exists\, l \in [k] \mid f(\bar X_l) \leq f(x^*) + \delta \right) 
    \;\;\geq\;\; 
    1 - q_{\mathrm{obj}}^k.
  \]
\end{restatable}

% Since $\psi(p,k)$ increases as $\mu$ decreases, an intuitive implication is that when $p$ is closer to a binary solution $x^*$, the likelihood of obtaining a sample with a comparable objective value increases. 

\blue{
The above bound provides a quantitative relationship among the optimality tolerance $\delta$, the objective sensitivity $L_f$, the expected distance $\mu$ from an optimal binary solution, and the sample size $k$. In particular, it implies two useful insights. First, for fixed $\delta$, $L_f$, and $k$, decreasing $\mu$ increases the probability of generating a sample with an objective value within $\delta$ of the optimum, 
which suggests the use of the binary penalty in \eqref{eq:saddle} to promote concentration toward binary solutions as the algorithm progresses. Second, the strength of the underlying relaxation remains important in GFORS, since, when a tighter relaxation yields a fractional solution closer to an optimal binary solution, the resulting smaller $\mu$ leads to a stronger sampling guarantee. Thus, the quality of the fractional solution obtained from the continuous optimization is particularly relevant. In this regard, when $Q$ is PSD and $\rho=0$, Theorem~\ref{thm:convex} shows that the algorithm converges to the global optimum of the underlying continuous relaxation. In practice, however, $\mu$ depends on the unknown optimizer $x^*$ and is therefore generally unavailable a priori. Consequently, the bound is best interpreted as a theoretical characterization of the factors governing sampling quality, rather than as an exact probability certificate that can be evaluated during optimization.
}

The next theorem provides the probability bound for the feasibility guarantee. Without loss of generality, we assume that only inequality constraints $Ax\geq b$ are present in \eqref{eq:bip}.

\begin{restatable}{theorem}{smpfeas}\label{thm:feasibility}
Given $p \in [0,1]^n$ and $k \in \mathbb N$ in Algorithm~\ref{alg:gsample},
\blue{define
  \begin{equation}\label{eq:q-obj}
      q_{\mathrm{feas}} = \exp\left(-2\left[\frac{\gamma_+^2(p)}{\eta^2} - \frac{\log m}{2}\right]_+\right).
  \end{equation}}Then, the following bound holds
\[
\Pr\left(\exists l \in [k] \mid A \bar X_{l} \geq b\right) \geq 1-q_{\mathrm{feas}}^k,
\]
where $A \in \mathbb R^{m\times n}$, $\gamma(p) = \min_{j \in [m]}\{\iprod{A_j, p} - b_j\}$ is the minimum slack, $\gamma_+(p):= \max\{\gamma(p), 0\}$, and $\eta = \max_{j \in [m]}\|A_j\|_2$ denotes the maximum row-wise Euclidean norm of $A$.
\end{restatable}

% One implication of this theorem is that a larger minimum slack in the fractional solution $p$ increases the likelihood of generating feasible binary solutions. On the other hand, when the feasible region $\{x \in [0,1]^n \mid Ax \ge b\}$ occupies only a small volume of the hypercube, the theorem indicates a weaker guarantee of feasibility. In particular, in the presence of equality constraints \eqref{eq:bip02}, the probability of obtaining feasible samples becomes extremely small. Therefore, we devote the next subsection to overcome this sampling pitfall.

One implication of this theorem is that a larger minimum slack in the fractional solution $p$ increases the likelihood of producing feasible binary samples. On the other hand, when the feasible region $\{x \in [0,1]^n \mid Ax \ge b\}$ occupies only a small portion of the hypercube, the guarantee of feasibility becomes weaker. 
\blue{
Unlike $q_{\mathrm{obj}}$, all quantities defining $q_{\mathrm{feas}}$ are directly computable for a given iterate $p$, making the resulting feasibility guarantee explicitly evaluable during the algorithm. The bound may nevertheless be conservative, primarily due to the use of the union bound in its derivation.
Combining the bounds in Theorems~\ref{thm:optprob} and \ref{thm:feasibility} through a union bound yields the following joint guarantee for inequality constrained problems.}

\blue{
\begin{corollary}
Let the assumptions of Theorems~\ref{thm:optprob} and~\ref{thm:feasibility} hold, and suppose that the problem contains only inequality constraints $Ax\ge b$. 
Then, for $k$ independent samples $\bar X_1,\ldots,\bar X_k$ generated by Algorithm~\ref{alg:gsample},
\[
\Pr\left(\exists\, l\in[k]:A\bar X_l\ge b,\quad f(\bar X_l)\le f(x^*)+\delta \right)
\ge 1 - \min\{1,q_{\mathrm{obj}}+q_{\mathrm{feas}}\}^{k}.
\]
\end{corollary}
}

% In particular, under equality constraints \eqref{eq:bip02}, the probability of obtaining feasible samples can be extremely small. To address this sampling pitfall, we will introduce two methods in the next subsection. 

\blue{
Under the equality constraints in~\eqref{eq:bip02}, we have $\gamma_+(p)=0$ and hence $q_{\mathrm{feas}}=1$, making the bound in Theorem~\ref{thm:feasibility} trivial. Intuitively, in this setting, generic Bernoulli sampling may have an extremely low probability of generating feasible samples. To mitigate this limitation, we introduce two mechanisms in the next subsection.
}

\subsection{Sampling with Equality Constraints}
We introduce two mechanisms to enhance the feasibility guarantees of sampling under equality constraints \eqref{eq:bip02}. Each mechanism has distinct strengths and is suited to different settings.

\subsubsection{Totally Unimodular (TU) Reformulation}
\label{sec:tureform}
A matrix is totally unimodular if every square submatrix has determinant in $\{-1, 0, 1\}$. Many binary integer programs adopt TU substructures \citep{wei2013total,rebman1974total}. In \eqref{eq:bip}, suppose a TU submatrix $B_{J}$ is identified with the row index set $J$, the following reformulation is exact and can remove all equality constraints associated with $J$.

\begin{restatable}{theorem}{totuni}
  In \eqref{eq:bip}, given row and column index sets $J$ and $I$ such that $B_J$ and $d_J$ are integral, $B_J$ is TU, and $B_{JI}$ is invertible, define $s:= B_{JI}^{-1}d_J$, $S:=-B_{JI}^{-1}B_{J\bar I}$, and let $\bar J$ and $\bar I$ denote the complement indices in the rows and columns of $B$. The following reformulation of \eqref{eq:bip} is exact,
$$
\min_{x_{\bar I} \in \{0,1\}^{n - |I|}} \left\{\iprod{x_{\bar I}, Q'x_{\bar I}} + \iprod{c', x_{\bar I}} + c'_0 \mid A' x_{\bar I} \geq b',\  B' x_{\bar I}  = d',\ Sx_{\bar I} \geq -s, \ Sx_{\bar I} \leq 1 - s\right\}
$$
where the model parameters $(Q',A',B',b',c',d', c'_0)$ are defined as
$$
\begin{aligned}
  Q' &:= S^\intercal Q_{II}  S + S^\intercal Q_{I\bar I} + Q_{I\bar I}^\intercal S + Q_{\bar I \bar I};\ A' := A_{\cdot I} S + A_{\cdot \bar I};\  B':= B_{\bar J I}S + B_{\bar J \bar I};\ b' := b - A_{I}s;\\
  c' &= 2 S^\intercal Q_{II} s + 2 Q^\intercal_{I\bar I}s + S^\intercal c_I + c_{\bar I};\ d' := d_{\bar J} - B_{\bar J I} s;\ c'_0 := \iprod{s, Q_{II}s} + \iprod{c_I, s}.
\end{aligned}
$$
\end{restatable}

%\begin{remark}
%  In this reformulation, the matrices $Q'$, $A'$, and $B'$ often lose sparsity because of the inverse in $S$, leading to a substantial memory burden. In the GPU implementation of the \textsc{TUReformulate} subroutine in Algorithm~\ref{algo:gfors}, we instead store the LU decomposition of $B_{JI}$ and use it to apply any multiplication of the form $B_{JI}^{-1}x$ efficiently during the gradient update steps.
%\end{remark}

This reformulation is particularly effective when the entire matrix $B$ is totally unimodular (see Section~\ref{sec:fl}). \blue{
Moreover, the sampling guarantees established in Section~\ref{sec:sanalysis} continue to apply after the TU reformulation and may yield tighter bounds by eliminating all equality constraints.
}
To address more general cases, we introduce a customized sampling method.

\subsubsection{Customized Sampling}
A customized sampling routine refers to any realization of \textsc{RandSampleStep} in Algorithm~\ref{algo:gfors} that, given the current fractional solution $p$ and a sample size $k$, returns the corresponding candidate solutions. An additional requirement is that this routine must be GPU-compatible. 
We illustrate this using the three-dimensional (3D) assignment problem.
\begin{subequations}
  \label{eq:assign3d}
\begin{align}
  \min_{x \in \{0,1\}^{n^3}} & \iprod{c, x} \\
\text{s.t.} 
                             & \sum_{(j,k) \in [n]^2} x_{ijk} = 1, \quad \forall i \in [n], \\
                             & \sum_{(i,k) \in [n]^2} x_{ijk} = 1, \quad \forall j \in [n], \\
                             & \sum_{(i,j) \in [n]^2} x_{ijk} = 1, \quad \forall k \in [n].
\end{align}
\end{subequations}
We can still apply the TU reformulation to eliminate the first two sets of constraints, but $n$ equality constraints remain. Algorithm~\ref{alg:3dassign} introduces a customized sampling subroutine that guarantees feasibility. The main idea is to first construct a feasible 3D assignment from $p$, and then apply parallelized local improvements using feasibility-preserving \emph{pairwise interchange} operations \citep{balas1991algorithm}.
\begin{algorithm}[t]
\caption{Customized Sampling for 3D Assignment}
\label{alg:3dassign}
\begin{algorithmic}
  \State Parameters: entry size multiple $\gamma$; improvement steps $L$
\Function{RandSampleStep}{$p, k$}
\State $p_{\gamma n} \gets $ \blue{largest $\min\{n^3, \lceil\gamma n\rceil\}$}  entries in $p$
\State $x_{\text{partial}} \gets$ partial non-conflict assignment according to the sorted $p_{\gamma n}$
\For{$l \in [k]$}
\State $\bar X_l \gets$ random completion of $x_{\text{partial}}$
\State $\bar X_l \gets$ local improvement of $\bar X_l$ with $L$ steps
\EndFor
  \State \Return $\bar X \in \{0,1\}^{k \times n}$
\EndFunction
\end{algorithmic}
\end{algorithm}

\blue{
In particular, for the partial non-conflict assignment, we retain the \(K=\min\{n^3,\lceil \gamma n\rceil\}\) largest entries of \(p\) and process them in decreasing order. At each step, we accept the current triple \((i,j,k)\) if none of its indices has been used, and then discard all remaining triples sharing \(i\), \(j\), or \(k\). This continues until \(n\) triples have been selected or the top-\(K\) list is exhausted. Any remaining indices are completed by randomly pairing the unused \(i\)-, \(j\)-, and \(k\)-indices without repetition. The top-\(K\) selection and ordering require \(O(n^3+K\log K)\) time, while the assignment-elimination scan requires \(O(K)\) time. For each feasible assignment, we then perform \(L\) local-improvement steps in parallel on the GPU. In each step, the selected triples are randomly paired, and for each pair, either the \(j\)-indices or the \(k\)-indices are swapped with equal probability. The swap is retained only if it reduces the total cost. The construction preserves feasibility: the partial assignment contains no repeated \(i\)-, \(j\)-, or \(k\)-indices, random completion uses each remaining index exactly once, and the pairwise interchange operations preserve these assignment constraints. In our implementation, we set \(\gamma=60\) and \(L=1000\).
}

% We construct a partial solution from the largest $\gamma n$ entries of $p$, which both reduces sorting cost and retains the most informative entries for candidate generation. This design enables a significantly more GPU-efficient implementation. The performance of this customized sampling subroutine will be presented in Section~\ref{sec:assign3d}.

A similar idea applies to other customized sampling designs: first generate a feasible candidate from the fractional solution $p$, then perform parallel local searches to explore additional candidates.
\blue{
However, due to the problem-specific nature of these subroutines, the sampling guarantees established in Section~\ref{sec:sanalysis} do not directly apply.
}

\begin{remark}[\blue{GFORS with Terminal Refinement}]
\blue{
Although the focus of this paper is the first-order and sampling framework, heuristic refinement procedures, as commonly used in methods such as B\&B, cuOpt, and CHAP, can be incorporated to further improve sampled binary solutions. The customized sampling scheme for 3D assignment in Algorithm~\ref{alg:3dassign} can be viewed as an iterative repair/refinement procedure applied after each sampling step. A simpler alternative is \emph{terminal refinement}, in which GFORS retains a set of top sampled candidates and applies a parallel refinement to them at the end.
}
\end{remark}

% \subsubsection{Monotone Relaxation}
% It is well-known that when the objective function $f$ is monotone, the binary solution space can be relaxed to its upper closure, i.e., all supersets of feasible solutions (for example, all $s$-$t$ connected subgraphs in the shortest path problem). Such relaxations often allow equality constraints to be replaced by inequalities. 
% 
% For instance, if $c$ is nonnegative or nonpositive in the 3D assignment formulation, all equalities can be relaxed to inequalities without changing the optimal solution (the same solution remains optimal under the relaxation). Our sampling algorithm can then operate on this relaxed formulation to generate candidate supersets of true assignment solutions, followed by a customized repair step to recover a valid assignment. 
% 
% This approach extends naturally to the case where $f$ is bimonotone, i.e., it becomes monotone after flipping certain entries of $x$.

\section{Implementation Details}
\label{sec:implement}

This section provides the implementation details of the remaining subroutines in Algorithm \ref{algo:gfors}. Apart from the straightforward \textsc{EvalBest}, which performs feasibility testing and optimality comparison via matrix multiplications, the others include: \textsc{Preprocess}, which normalizes the input parameters to handle heterogeneous problem instances; \textsc{UpdatePenalty}, which gradually increases the binary gap penalty multiplier $\rho$; and \textsc{CheckHalt}, which enforces the stopping criteria.

% \paragraph{\textsc{Preprocess.}} This subroutine performs normalization on the model parameters $(Q, A, B, b, c, d, c_0)$, potentially the ones obtained after applying the \textsc{TUReformulate} subroutine. This is a common step used in various mathematical programming solvers to ensure a consistent algorithm performance across distinct instances \cite{gould2004preprocessing,koch2022progress,wolsey2020integer}. 

\paragraph{\textsc{Preprocess.}}
This subroutine normalizes the model parameters $(Q, K, c, r)$, which may include those produced by the \textsc{TUReformulate} subroutine. Such normalization is a standard procedure in many mathematical programming solvers to ensure consistent algorithmic performance across distinct instances \citep{gould2004preprocessing,koch2022progress}. Specifically, the following steps are performed in this function:
\begin{itemize}
  \item Normalize each row in $K$ by its 2-norm and update $r$ accordingly.
  \item Normalize $Q$ and $c$ by $\|Q\|_2 + \|c\|_2$ (i.e., spectral norm of $Q$ plus 2-norm of $c$).
  \item Normalize $K$ and $r$ by the spectral norm of $K$.
\end{itemize}
\blue{
The purpose of the last step is to normalize away the dependence on $\|K\|$ in the step-size product, allowing the practical tuning rule to be expressed as $\tau_1\tau_2 \leq \sigma$.
  For efficient spectral-norm estimation of large sparse matrices, we use \texttt{scipy.sparse.linalg.svds} to approximate the largest singular value without forming a dense matrix.
}
In our implementation, we use the tuning parameter $\sigma \in (0,1)$ to control the exploration aggressiveness via $\tau_1 \tau_2 \leq \sigma$. 
\blue{To maximize the primal and dual step sizes under this condition, we initialize $\tau_1=\tau_2=\sqrt{\sigma}$ and heuristically update
\[
r_k := \operatorname{clip}\!\left(\frac{f_k}{d_k},\,0.1,\,10\right), 
\qquad
\tau_1^k \gets \sqrt{r_k\sigma},
\qquad
\tau_2^k \gets \sqrt{\frac{\sigma}{r_k}},
\]
where $f_k$ and $d_k$ denote the primal and dual residuals at iteration $k$, respectively.}

\paragraph{\textsc{UpdatePenalty.}} 
The strength of the binary penalty parameter $\rho$ dictates how strongly the current iterate $x_k$ is attracted toward binary solutions. 
Intuitively, we initialize $\rho$ with a small value so that the iterate trajectory first follows the relaxation of \eqref{eq:bip}, allowing $x_k$ to move toward a global optimum when $Q$ is positive semidefinite. 
Subsequently, $\rho$ is gradually increased to encourage convergence to binary solutions in later stages, thereby boosting the probability of obtaining a high-quality feasible solution as suggested by Theorem~\ref{thm:optprob}. 
With tunable parameters $(\rho_{\min}, \rho_{\max}, T, p, \delta)$, we employ the heuristic update 
\begin{equation}
\label{eq:pscheduler}
\blue{
    \tilde \rho_n  = \rho_{\min}\left(1 + \tfrac{n}{T}\right)^p, \qquad\rho_{n} = \text{clip}(\tilde \rho_n, \rho_{n-1} + \delta, \rho_{\max}),}
\end{equation}
where $\rho_{\min}$ and $\rho_{\max}$ set the lower and upper bounds of the penalty, $T$ and $p$ control the growth rate, and $\delta$ enforces a minimum increment pace. 
Here, $n$ is a counter that increases only after designated iterations. 
This update strategy balances early-stage exploration guided by the continuous relaxation with late-stage exploitation that strengthens the pull toward integrality, thereby balancing global search and binary feasibility.

\paragraph{\textsc{CheckHalt.}} 
Since solution generation primarily relies on the sampling subroutine, the stopping criterion is designed to detect whether further improvement is unlikely, thereby enabling potential early termination. 
With the aid of a value-stalling detector, this subroutine returns \texttt{true} once all three indicators, i.e., the primal feasibility gap, the dual feasibility gap, and the binary gap, are either satisfied within tolerance or have remained stalled for a sufficient number of iterations.

\paragraph{\textsc{TUReformulate.}} \blue{
In our implementation, the row and column index sets \(J\) and \(I\) used for the TU reformulation are specified by the user. Given \(J\) and \(I\), we partition the selected equality constraints as \(B_{JI}x_I + B_{J\bar I}x_{\bar I} = d_J\) and apply the substitution \(x_I = B_{JI}^{-1}(d_J - B_{J\bar I}x_{\bar I})\) implicitly, without explicitly forming  \(B_{JI}^{-1}\).
In particular,
after the normalization in the preprocessing step, \(B_{JI}\) is factorized once on the GPU using PyTorch's pivoted LU factorization routine \texttt{torch.linalg.lu\_factor}. Reconstructions of \(x_I\), including batched candidate reconstruction, use \texttt{torch.linalg.lu\_solve} with the stored factors. Numerical stability is handled through the preprocessing normalization and the two routines provided by PyTorch.
}

% \paragraph{\textsc{CheckHalt.}} Since our solution generation mainly comes from the sampling subroutine, the stopping criteria is based on whether obvious improvement is possible for potential early stop. With the help of a value stalling detector, this subroutine returns true whenever all three indicators---primal feasibility gap, dual feasibility gap, and binary gap---are either closed within tolerance or stalled for sufficient iterations.

\section{Numerical Experiments}
\label{sec:expr}
In this section, we evaluate GFORS on a collection of classic binary integer programs and benchmark its performance against the latest versions of Gurobi \blue{ and cuOpt} available at the time of experimentation. We begin with an aggregated view of all computational results, followed by problem-specific analyses. Because integer programs often exhibit heterogeneous performance across problem classes, these detailed results offer additional insight into GFORS and its individual components introduced in Sections~\ref{sec:framework} and~\ref{sec:implement}.
\blue{To ensure reproducibility, all datasets are publicly available on Hugging Face\footnote{Instance repository \url{https://huggingface.co/datasets/NingjiWei/GFORS_BIP_Instances}}, while the complete codebase and detailed experiment settings are available on GitHub\footnote{Codebase repository \url{https://github.com/tidues/GFORS}}.}

\paragraph{Experiment Setup.} 
All experiments were conducted on a Slurm-managed HPC cluster (Slurm 24.05) running Red Hat Enterprise Linux 9.4 with kernel 5.14.0. Each job was allocated 16 CPU cores from an Intel Xeon Gold 6348Y processor, one NVIDIA L40S GPU (48 GB, CUDA 12.5, driver 555.42.02), and 96 GB RAM, matching the standard memory ratio used in the literature \citep{lu2023cupdlp}. We developed GFORS in Python 3.12.11 with PyTorch 2.5.1 for GPU acceleration, and used Gurobi \blue{13.0.2 and cuOpt 26.08} as the two baseline BIP solvers. 
\blue{Unless stated otherwise, GFORS initializes with \(x^0=0.5 \times \mathbf{1}\) and all dual variables at $\mathbf{0}$, performs $k_r=10$ sampling rounds of $k_b=10$ Bernoulli candidates ($100$ samples) every $k_{\mathrm{int}}=50$ iterations, and uses \(\sigma=0.3\).}
Each instance was run with a time limit of $1{,}800$ seconds. 
\blue{The Gurobi and cuOpt baselines use the same time limit, with all other parameters left at their default values. Reported solver times exclude model construction, preprocessing, and post-processing, since these implementation-dependent operations are not intrinsic to the optimization algorithm and may be shared, reused, or handled differently across solvers. Whenever applicable, TU-reformulation time is reported separately. 
For the penalty scheduler~\eqref{eq:pscheduler}, we set
$(\rho_{\min}, \rho_{\max}, T, p, \delta) = (0.01, 100, 1000, 1.01, 0.01)$.}

% \blue{The Gurobi and cuOpt baselines use the same time limit, with all other parameters remaining at their default values. Model-construction, preprocessing, and post times are excluded from the reported solver times, since ... Whenever applicable, TU-reformulation time is reported separately.}

\paragraph{Problem Instances.} 
We evaluate our method on five classic BIP classes: set cover, knapsack, max-cut, 3D assignment, and facility location.
For each class, we scale the instance size by powers of two until the number of nonzeros (NNZ) in the model parameters exceeds $2 \times 10^9$. 
\blue{We treat the quadratic and linear integer max-cut models as two formulations of the same graph instances. The linear formulation has one fewer size because its largest model exceeds the $2\times10^9$ NNZ limit.}
At each size, we generate five \blue{independent random} instances, yielding \blue{$455$} test instances in total. 
\blue{Each solver configuration is run once on each instance, so the reported averages are taken over the five independently generated instances at each size.}
This exponential scaling systematically evaluates the scalability of the GFORS framework. 
We also include an additional experiment on the traveling salesman problem \blue{(TSP) with $45$ instances} as a stress test in Appendix~\ref{sec:tsp}.

\subsection{Aggregated Results Overview}
\label{sec:agg}
Since the instances span multiple problem classes, we use the total NNZ in the model parameters $(Q, K)$ to define instance size, and aggregate computational performance into 28 groups based on $\lfloor \log_{2}(\mathrm{nnz}) \rfloor$. Because the problem classes differ in structure, the number of instances per group is not perfectly uniform. With the exception of the smallest group ($2^3$) containing only five instances, each NNZ group includes between 10 and 20 instances drawn from multiple classes. 

% To ensure fair comparison, we exclude the TSP class from the aggregated results, as it differs fundamentally from the other five problem classes (lacking a natural compact BIP formulation), and GFORS solves only the smallest TSP instances, which would otherwise skew the aggregation result. A detailed analysis of the TSP class is provided in Section~\ref{sec:tsp}.

Though our GFORS implementation allows for many tunable parameters, in this experiment we adjust only four settings: the step-size parameter $\sigma \in (0,1)$, which controls exploration aggressiveness; the sampling batch size $k_b$, which can be used to limit peak GPU memory usage; the TU reformulation, applied when row and column index sets $J$ and $I$ are given (see Section~\ref{sec:tureform}); and customized sampling subroutine if provided.

We benchmark GFORS using a fixed set of configurations against the default settings of Gurobi with presolve enabled and disabled, \blue{as well as cuOpt with presolve enabled}. Because GFORS is not an exact solver but is designed to identify high-quality solutions for large instances, we evaluate performance using two complementary metrics: (i) \emph{solved percentage}, the fraction of instances in each NNZ group for which a feasible solution is obtained within the time limit; and (ii) \emph{time-to-target}, which provides a fair basis of comparison since Gurobi \blue{ and cuOpt} may continue refining incumbents long after producing strong early solutions. For GFORS, time-to-target is the time to its best solution before halting; for Gurobi \blue{ and cuOpt}, it is the \blue{first time a solution reaches or surpasses} the GFORS incumbent \blue{($\leq$ for minimization, $\geq$ for maximization)}, or the time limit if the target is not reached. Metric (i) captures each algorithm’s ability to find feasible solutions, while (ii) measures the time required for \blue{the two baselines} to reach GFORS's incumbent, thereby indicating whether GFORS can effectively fulfill its purpose of producing high-quality feasible solutions, especially on large-scale instances.

\begin{figure}[t]
  \centering
  \begin{subfigure}{0.485\textwidth}
    \centering
    \includegraphics[width=\linewidth]{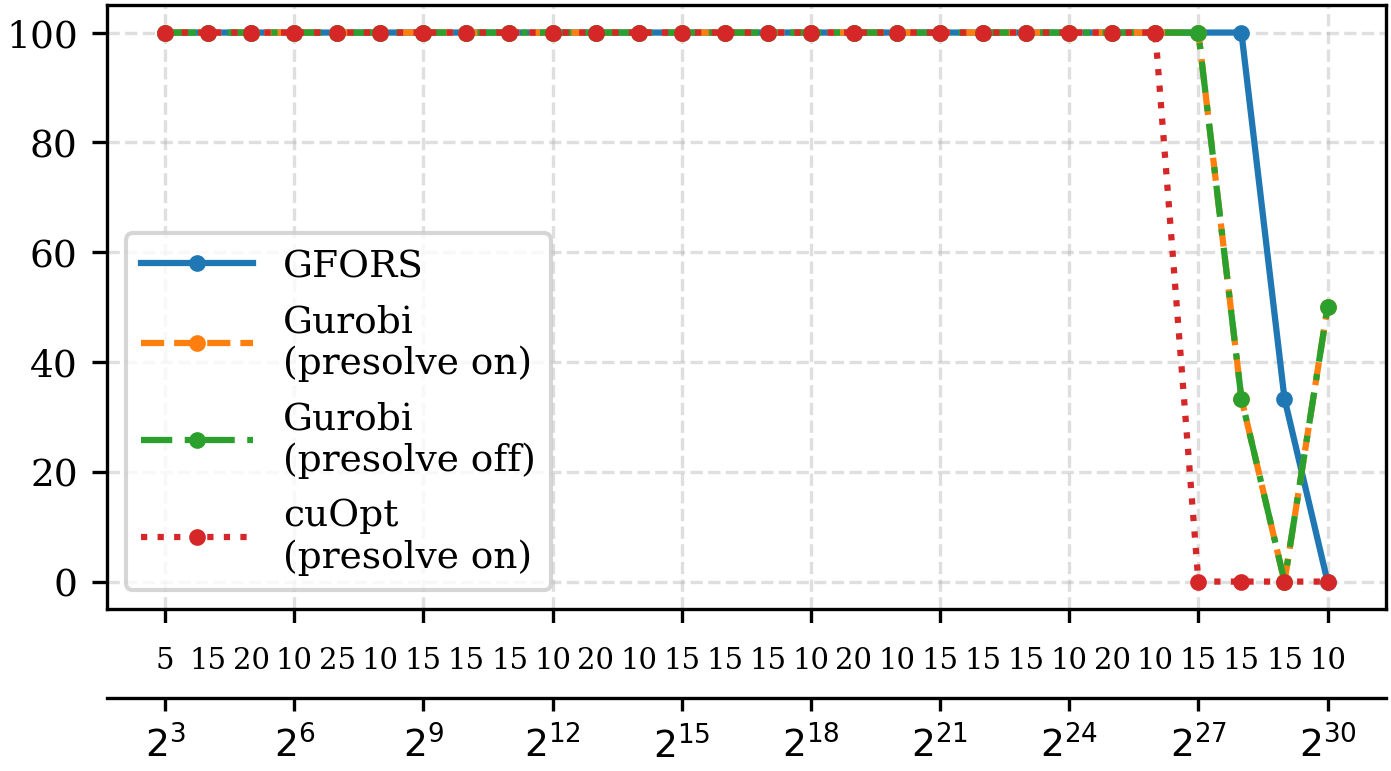}
    \caption{Solved percentage (\%) vs NNZ}
    \label{fig:perf-a}
  \end{subfigure}\hfill
  \begin{subfigure}{0.495\textwidth}
    \centering
    \includegraphics[width=\linewidth]{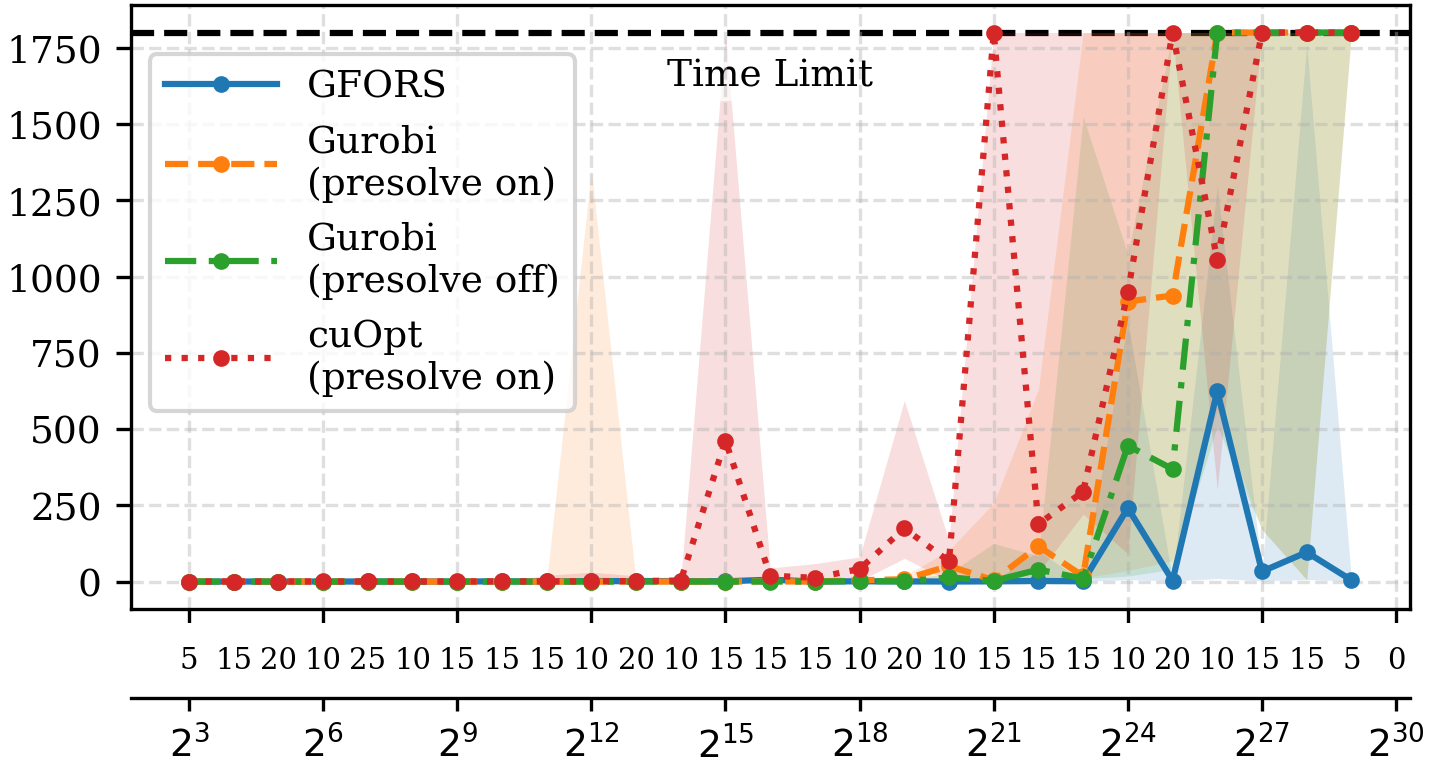}
    \caption{Time-to-target (sec.) vs NNZ}
    \label{fig:perf-b}
  \end{subfigure}
\caption{
Performance vs.\ NNZ. Solved percentages and time-to-target values are aggregated by $\lfloor \log_{2}(\mathrm{nnz}) \rfloor$\blue{; the label below each data point gives the number of aggregated instances.}
Because Gurobi \blue{and cuOpt} may spend substantial time refining incumbents after finding strong early solutions, we use time-to-target as the comparison metric. For GFORS, this is the time to its best solution before termination; for Gurobi \blue{and cuOpt}, it is the first time a solution \blue{reaches or surpasses the GFORS objective ($\leq$ for minimization, $\geq$ for maximization)}, or $1{,}800$ seconds if no such solution is found. \blue{Terminal refinement is not used for the GFORS results in this figure.}}

% \caption{Performance vs.\ number of nonzeros (NNZ). Solved percentages and times are aggregated by $\lfloor \log_{2}(\mathrm{nnz}) \rfloor$. \emph{Times-to-target} are shown as medians with shaded interquartile ranges: for GFORS, this is the time to its best solution before halting; for Gurobi (w/ and w/o presolve), it is the first time a solution improves upon the GFORS optimum, or $1{,}800$ seconds if none is found.}
  \label{fig:perf-both}
\end{figure}

The aggregated results are presented in Figure~\ref{fig:perf-both}, where the $x$-axis shows NNZ on a logarithmic scale, \blue{with each data point labeled by the number of instances included in the aggregation}. The two subfigures report the solved percentage and time-to-target, respectively.

Figure~\ref{fig:perf-a} shows that GFORS maintains a $100\%$ solved rate across a wider range of NNZ groups than the baselines, with degradation appearing only at the largest scales. \blue{Both Gurobi configurations exhibit a similar sharp decline at large NNZ, while cuOpt begins to fail earlier and solves no instances beyond the scale of $2^{27}$. 
The rebound of the two Gurobi variants at the largest scale is driven entirely by the set cover instances (Table~\ref{tab:sc}), where both obtain high-quality root-node incumbents but make no further improvement.
}

Figure~\ref{fig:perf-b} reports only those problem instances for which GFORS obtained a feasible solution. The plot shows median time-to-target values as line curves, with shaded bands indicating the first and third quartiles. \blue{Both Gurobi variants generally outperform GFORS on small- to medium-sized instances, while cuOpt shows greater variability and occasional large time-to-target values even at intermediate scales.} In the range $2^{11}$ to $2^{15}$, \blue{the baselines exhibit noticeable time-to-target variation on certain instances,} reflected in the large third-quartile values. For large-scale instances, however, GFORS consistently finds high-quality solutions quickly, whereas \blue{both Gurobi variants and cuOpt take longer to reach the same targets and frequently hit the time limit.} These results highlight GFORS’s efficiency on large-scale instances relative to \blue{the baselines.}

\subsection{Problem-Specific Analysis}

\blue{This section analyzes results by problem class, using Gurobi and cuOpt with presolve enabled as baselines, since presolve generally matches or improves performance over the presolve-off variants.}
Instance inputs are split into (i) size parameters that determine model scale and (ii) randomization parameters (e.g., edge weights, item costs). For each size configuration, tables report the average best objective and the associated average runtime over five random instances. 
% To isolate core solver performance, we exclude model-construction, preprocessing, and postprocessing times for all methods (e.g., \textsc{Preprocess} in Algorithm~\ref{algo:gfors}), since parsing/instantiation costs are interface-dependent, orthogonal to algorithmic behavior, and typically small relative to runtime. TU reformulation time is reported separately wherever applicable.

\subsubsection{Set Cover: Effect of the Step-Size Parameter $\sigma$}

\begin{table}[!tbp]
  \centering
  \blue{
 \scriptsize
\setlength{\tabcolsep}{4pt}
\begin{tabular}{lrrrrrrrrrrr}
\toprule
Config & \multirow{2}{*}{NNZ}
& \multicolumn{5}{c}{Avg. Objective Value}
& \multicolumn{5}{c}{Avg. Runtime (sec.)} \\
\cmidrule(lr){3-7}\cmidrule(lr){8-12}
$(m,n)$ & & GRB & cuOpt & $\sigma=0.99$ & $\sigma=0.5$ & $\sigma=0.3$
& GRB & cuOpt & $\sigma=0.99$ & $\sigma=0.5$ & $\sigma=0.3$ \\
\midrule
$(5,10)$ & $1.74e1$ & $7.8$ & $7.8$ & $7.8$ & $7.8$ & $7.8$ & $0.00$ & $0.35$ & $0.74$ & $0.93$ & $0.78$ \\
$(10,20)$ & $3.30e1$ & $15.6$ & $15.6$ & $15.6$ & $15.6$ & $15.6$ & $0.00$ & $0.35$ & $0.78$ & $1.09$ & $1.21$ \\
$(20,40)$ & $6.82e1$ & $35.0$ & $35.0$ & $35.2$ & $35.2$ & $35.2$ & $0.00$ & $0.35$ & $1.93$ & $2.96$ & $3.19$ \\
$(40,80)$ & $1.40e2$ & $57.0$ & $57.0$ & $57.2$ & $57.2$ & $57.2$ & $0.00$ & $0.31$ & $1.77$ & $1.73$ & $1.25$ \\
$(80,160)$ & $4.55e2$ & $73.2$ & $73.2$ & $73.8$ & $73.8$ & $74.2$ & $0.01$ & $0.32$ & $1.85$ & $3.39$ & $5.69$ \\
$(160,320)$ & $1.82e3$ & $78.8$ & $78.8$ & $81.4$ & $83.4$ & $84.0$ & $0.03$ & $0.48$ & $4.09$ & $5.01$ & $6.72$ \\
$(320,640)$ & $7.04e3$ & $77.4$ & $77.4$ & $86.8$ & $89.2$ & $91.8$ & $0.90$ & $1.14$ & $4.54$ & $5.70$ & $6.96$ \\
$(640,1280)$ & $2.86e4$ & $72.4$ & $72.4$ & $98.4$ & $100.8$ & $107.8$ & $258.53$ & $671.38$ & $5.51$ & $7.17$ & $8.16$ \\
$(1280,2560)$ & $1.15e5$ & $73.9$ & $72.8$ & $170.0$ & $199.8$ & $219.4$ & $1800.07$ & $1800.52$ & $8.26$ & $10.63$ & $12.50$ \\
$(2560,5120)$ & $4.63e5$ & $82.5$ & $82.0$ & $355.2$ & $291.8$ & $389.8$ & $1800.17$ & $1800.59$ & $15.53$ & $4.32$ & $1.27$ \\
$(5120,10240)$ & $1.84e6$ & $100.0$ & $101.4$ & $463.8$ & $953.8$ & $377.4$ & $1803.95$ & $1801.58$ & $1.28$ & $1.29$ & $3.28$ \\
$(10240,20480)$ & $7.36e6$ & $122.2$ & $124.6$ & $397.8$ & $370.8$ & $241.2$ & $1800.91$ & $1801.48$ & $1.38$ & $43.61$ & $315.06$ \\
$(20480,40960)$ & $2.94e7$ & $604.1$ & $142.4$ & $322.0$ & $278.0$ & $273.4$ & $1801.47$ & $1801.02$ & $148.84$ & $707.61$ & $1262.19$ \\
$(40960,81920)$ & $1.17e8$ & $1120.6$ & $165.4$ & $314.4$ & $953.2$ & $1615.4$ & $1806.47$ & $1802.85$ & $1625.92$ & $1800.02$ & $1800.01$ \\
$(81920,163840)$ & $4.69e8$ & $1271.2$ & $\mathrm{OOM}$ & $3689.2$ & $8825.0$ & $15084.4$ & $992.47$ & $52.78$ & $1800.05$ & $1800.02$ & $1800.04$ \\
$(163840,327680)$ & $1.88e9$ & $1427.6$ & $\mathrm{OOM}$ & $\mathrm{OOM}$ & $\mathrm{OOM}$ & $\mathrm{OOM}$ & $162.35$ & $10.23$ & -- & -- & -- \\
\bottomrule
\end{tabular}

}
\caption{Set Cover: objective values and runtimes. The parameter $\sigma$ controls the first-order step-size schedule in GFORS; \blue{$\mathrm{TL}$ and $\mathrm{OOM}$ indicate no feasible solution before the time or memory limit, respectively.}}
  \label{tab:sc}
\end{table}

The set cover problem seeks the minimum-cost collection of subsets whose union covers the entire ground set, formulated as
\[
  \min_{x \in \{0,1\}^n} \left\{ \iprod{c, x} \,\middle|\, \sum_{i \in S_j} x_i \geq 1,\ \forall j \in [m]\right\}.
\]
The test instances range from $(m,n)=(5,10)$ up to $(81{,}920,163{,}840)$, doubling both $m$ and $n$ at each step. \blue{For each covering row, the number of selected columns is drawn uniformly between approximately $0.02n$ and $0.05n$, with a minimum row cardinality of three, and the corresponding columns are sampled uniformly without replacement. All selected coefficients are set to one, while the item costs are sampled independently and uniformly from $\{1,\ldots,10\}$.} We evaluate three GFORS configurations with $\sigma \in \{0.3,0.5,0.99\}$, where $\sigma$ controls the step-size aggressiveness through $\tau_1 \tau_2 \leq \sigma$. Detailed results are reported in Table~\ref{tab:sc}.

% The test instances range from $(m,n) = (5,10)$ up to $(81{,}920,163{,}840)$, doubling both $m$ and $n$ at each step. The cost vector $c$ and sets $S_j$ are randomly generated. We evaluate three GFORS configurations with $\sigma \in \{0.3,0.5,0.99\}$, where $\sigma$ controls the step-size aggressiveness through $\tau_1 \tau_2 \leq \sigma$. Detailed results are reported in Table~\ref{tab:sc}. 

\blue{Compared with Gurobi, GFORS trades some solution quality for shorter runtimes on small- to medium-sized instances. From approximately $\mathrm{NNZ}=10^5$ to $10^7$, GRB obtains better objectives but typically reaches the time limit, whereas GFORS terminates within seconds to minutes. At $\mathrm{NNZ}=2.94\times10^7$ and $1.17\times10^8$, GFORS also achieves better objectives than GRB. At the two largest configurations, GRB again obtains better objectives because its presolve identifies strong root-node incumbents. cuOpt closely tracks GRB on small and medium instances and outperforms the other methods at $m=20{,}480$ and $40{,}960$, but fails to find a feasible solution at the two largest sizes attempted.}
Across GFORS settings, on medium- to large-sized instances, smaller $\sigma$ values (e.g., $0.3,0.5$) generally achieve better objectives, whereas on the largest instances, larger $\sigma$ (e.g., $0.99$) is preferable due to its greater search efficiency (large step size) under the time limit.

\paragraph{Effects of Terminal Refinement.}

\blue{
To evaluate the effect of combining GFORS with terminal refinement, we retain up to the top $500$ sampled candidates from each GFORS run and apply a parallel GPU refinement procedure afterward. Specifically, each selected set is examined once in descending order of cost and removed whenever feasibility is preserved. 
Table~\ref{tab:screpair} reports the results, showing that terminal refinement has little effect on the smallest instances but becomes increasingly beneficial as problem size grows. From approximately $\mathrm{NNZ}=10^6$ onward, significant objective improvements are obtained across the three $\sigma$ settings, reaching over $90\%$ at the largest solvable size. The additional refinement cost remains modest relative to the GFORS solve, requiring at most about $19$ seconds. This indicates that combining terminal refinement with GFORS can strengthen incumbent quality at large scales with limited additional computational effort.
}

\begin{table}[!tbp]
  \centering
  \blue{
 \scriptsize
\begin{tabular}{lrrrrrrr}
\toprule
Config & \multirow{2}{*}{NNZ}
& \multicolumn{3}{c}{Avg. Repaired Objective Value (Impr.)}
& \multicolumn{3}{c}{Avg. Repair Time (sec.)} \\
\cmidrule(lr){3-5}\cmidrule(lr){6-8}
$(m,n)$ & 
& $\sigma=0.99$ & $\sigma=0.5$ & $\sigma=0.3$
& $\sigma=0.99$ & $\sigma=0.5$ & $\sigma=0.3$ \\
\midrule
$(5,10)$
& $1.74e1$
& $7.8~(-0.00\%)$ & $7.8~(-0.00\%)$ & $7.8~(-0.00\%)$
& $0.003$ & $0.003$ & $0.017$ \\

$(10,20)$
& $3.30e1$
& $15.6~(-0.00\%)$ & $15.6~(-0.00\%)$ & $15.6~(-0.00\%)$
& $0.005$ & $0.006$ & $0.025$ \\

$(20,40)$
& $6.82e1$
& $35.2~(-0.00\%)$ & $35.0~(-0.57\%)$ & $35.2~(-0.00\%)$
& $0.010$ & $0.010$ & $0.032$ \\

$(40,80)$
& $1.40e2$
& $57.2~(-0.00\%)$ & $57.2~(-0.00\%)$ & $57.2~(-0.00\%)$
& $0.014$ & $0.015$ & $0.037$ \\

$(80,160)$
& $4.55e2$
& $73.6~(-0.27\%)$ & $73.6~(-0.27\%)$ & $73.6~(-0.81\%)$
& $0.020$ & $0.023$ & $0.045$ \\

$(160,320)$
& $1.82e3$
& $80.2~(-1.47\%)$ & $80.2~(-3.84\%)$ & $81.4~(-3.10\%)$
& $0.028$ & $0.031$ & $0.056$ \\

$(320,640)$
& $7.04e3$
& $83.6~(-3.69\%)$ & $84.6~(-5.16\%)$ & $85.6~(-6.75\%)$
& $0.038$ & $0.040$ & $0.054$ \\

$(640,1280)$
& $2.86e4$
& $93.0~(-5.49\%)$ & $97.2~(-3.57\%)$ & $106.2~(-1.48\%)$
& $0.044$ & $0.042$ & $0.052$ \\

$(1280,2560)$
& $1.15e5$
& $163.8~(-3.65\%)$ & $191.6~(-4.10\%)$ & $210.4~(-4.10\%)$
& $0.042$ & $0.043$ & $0.054$ \\

$(2560,5120)$
& $4.63e5$
& $339.6~(-4.39\%)$ & $147.4~(-49.49\%)$ & $144.2~(-63.01\%)$
& $0.045$ & $0.142$ & $0.242$ \\

$(5120,10240)$
& $1.84e6$
& $130.8~(-71.80\%)$ & $129.2~(-86.45\%)$ & $128.4~(-65.98\%)$
& $0.325$ & $0.515$ & $0.889$ \\

$(10240,20480)$
& $7.36e6$
& $151.0~(-62.04\%)$ & $150.8~(-59.33\%)$ & $149.2~(-38.14\%)$
& $1.038$ & $1.574$ & $0.395$ \\

$(20480,40960)$
& $2.94e7$
& $174.6~(-45.78\%)$ & $173.4~(-37.63\%)$ & $173.6~(-36.50\%)$
& $3.191$ & $0.267$ & $0.438$ \\

$(40960,81920)$
& $1.17e8$
& $198.4~(-36.90\%)$ & $198.8~(-79.14\%)$ & $198.6~(-87.71\%)$
& $0.413$ & $1.437$ & $2.242$ \\

$(81920,163840)$
& $4.69e8$
& $223.4~(-93.94\%)$ & $224.6~(-97.45\%)$ & $224.8~(-98.51\%)$
& $7.037$ & $13.054$ & $18.928$ \\

$(163840,327680)$
& $1.88e9$
& $\mathrm{OOM}$ & $\mathrm{OOM}$ & $\mathrm{OOM}$
& -- & -- & -- \\
\bottomrule
\end{tabular}

}
\caption{Set Cover: GFORS with parallel terminal refinement applied to the top $500$ sampled candidates. \emph{Impr.} denotes the improvement relative to the GFORS-only objective in Table~\ref{tab:sc}.}
  \label{tab:screpair}
\end{table}

\subsubsection{Knapsack: Memory Control via Batch Size}

The knapsack problem searches for the subset of items with maximum total value subject to a weight capacity constraint, formulated as
\[
  \max_{x \in \{0,1\}^n} \left\{ \iprod{v, x} \,\middle|\, \iprod{w, x} \leq W \right\}.
\]
In our experiments, the number of items $n$ scales from $10$ up to $1.34 \times 10^9$, doubling at each step. \blue{Item values $v$ and weights $w$ are sampled independently and uniformly from the integers $\{1,\ldots,10\}$,} with $W$ fixed at half of the total weight. Detailed results are reported in Table~\ref{tab:knapsack}.

\begin{table}[!tbp]
  \centering
  \blue{
 \scriptsize
\begin{tabular}{rrrrrrrrr}
\toprule
Config & \multirow{2}{*}{NNZ}
& \multicolumn{3}{c}{Avg. Objective Value}
& \multicolumn{3}{c}{Avg. Runtime (sec.)}
& \multirow{2}{*}{GPU Mem (MB)} \\
\cmidrule(lr){3-5}\cmidrule(lr){6-8}
$n$ & & GRB & cuOpt & GFORS & GRB & cuOpt & GFORS & \\
\midrule
$10$ & $1.00e1$ & $3.80e1$ & $3.80e1$ & $3.51e1$ & $0.00$ & $0.38$ & $1.36$ & $30.17$ \\
$20$ & $2.00e1$ & $9.02e1$ & $9.02e1$ & $8.91e1$ & $0.00$ & $0.40$ & $1.24$ & $30.17$ \\
$40$ & $4.00e1$ & $1.67e2$ & $1.67e2$ & $1.66e2$ & $0.00$ & $0.42$ & $1.31$ & $30.17$ \\
$80$ & $8.00e1$ & $3.37e2$ & $3.37e2$ & $3.37e2$ & $0.00$ & $0.44$ & $1.28$ & $30.17$ \\
$160$ & $1.60e2$ & $6.89e2$ & $6.89e2$ & $6.88e2$ & $0.00$ & $0.41$ & $1.16$ & $30.17$ \\
$320$ & $3.20e2$ & $1.37e3$ & $1.37e3$ & $1.37e3$ & $0.00$ & $0.74$ & $0.79$ & $30.18$ \\
$640$ & $6.40e2$ & $2.79e3$ & $2.79e3$ & $2.79e3$ & $0.00$ & $0.97$ & $0.75$ & $30.20$ \\
$1,280$ & $1.28e3$ & $5.48e3$ & $5.48e3$ & $5.48e3$ & $0.00$ & $4.31$ & $0.84$ & $30.23$ \\
$2,560$ & $2.56e3$ & $1.10e4$ & $1.10e4$ & $1.10e4$ & $0.01$ & $42.12$ & $0.74$ & $30.31$ \\
$5,120$ & $5.12e3$ & $2.20e4$ & $2.20e4$ & $2.20e4$ & $0.01$ & $137.80$ & $1.28$ & $30.45$ \\
$10,240$ & $1.02e4$ & $4.41e4$ & $4.41e4$ & $4.41e4$ & $0.02$ & $1553.92$ & $0.89$ & $30.75$ \\
$20,480$ & $2.05e4$ & $8.81e4$ & $8.81e4$ & $8.80e4$ & $0.03$ & $3.61$ & $0.69$ & $31.33$ \\
$40,960$ & $4.10e4$ & $1.77e5$ & $1.77e5$ & $1.77e5$ & $0.06$ & $9.09$ & $0.99$ & $34.51$ \\
$81,920$ & $8.19e4$ & $3.54e5$ & $2.90e5$ & $3.54e5$ & $0.16$ & $21.19$ & $0.98$ & $40.87$ \\
$163,840$ & $1.64e5$ & $7.08e5$ & $4.69e5$ & $7.08e5$ & $0.28$ & $59.32$ & $0.86$ & $51.54$ \\
$327,680$ & $3.28e5$ & $1.42e6$ & $1.42e6$ & $1.42e6$ & $0.58$ & $83.15$ & $1.44$ & $69.25$ \\
$655,360$ & $6.55e5$ & $2.83e6$ & $2.83e6$ & $2.83e6$ & $0.99$ & $158.53$ & $1.15$ & $120.00$ \\
$1,310,720$ & $1.31e6$ & $5.66e6$ & $5.66e6$ & $5.66e6$ & $1.96$ & $298.25$ & $0.99$ & $230.16$ \\
$2,621,440$ & $2.62e6$ & $1.13e7$ & $1.13e7$ & $1.13e7$ & $3.99$ & $1356.34$ & $1.20$ & $404.83$ \\
$5,242,880$ & $5.24e6$ & $2.26e7$ & $2.26e7$ & $2.26e7$ & $7.76$ & $288.75$ & $1.76$ & $787.50$ \\
$10,485,760$ & $1.05e7$ & $4.53e7$ & $4.53e7$ & $4.53e7$ & $15.95$ & $1309.02$ & $3.00$ & $1547.50$ \\
$20,971,520$ & $2.10e7$ & $9.06e7$ & $4.40e3$ & $9.06e7$ & $31.41$ & $194.50$ & $8.88$ & $3112.83$ \\
$41,943,040$ & $4.19e7$ & $1.81e8$ & $4.41e3$ & $1.81e8$ & $63.99$ & $382.73$ & $20.49$ & $6200.83$ \\
$83,886,080$ & $8.39e7$ & $3.62e8$ & $4.00e1$ & $3.56e8$ & $130.90$ & $31.39$ & $25.91$ & $12302.16$ \\
$167,772,160$ & $1.68e8$ & $5.93e8$ & $\mathrm{OOM}$ & $7.25e8$ & $216.64$ & $32.37$ & $69.50$ & $24718.16$ \\
$335,544,320$ & $3.36e8$ & $\mathrm{OOM}$ & $\mathrm{OOM}$ & $1.45e9$ & -- & $25.73$ & $135.66$ & $49080.83$ \\
$671,088,640$ & $6.71e8$ & $\mathrm{OOM}$ & $\mathrm{OOM}$ & $\mathrm{OOM}$ & -- & $43.94$ & -- & -- \\
$1,342,177,280$ & $1.34e9$ & $\mathrm{OOM}$ & $\mathrm{OOM}$ & $\mathrm{OOM}$ & -- & $46.57$ & -- & -- \\
\bottomrule
\end{tabular}

  }
 \caption{Knapsack: objective values and runtimes. \emph{GPU Mem} reports peak GPU memory usage during computation; \blue{$\mathrm{TL}$ and $\mathrm{OOM}$ indicate no feasible solution before the time or memory limit, respectively.}}
  \label{tab:knapsack}
\end{table}

% In our experiments, the number of items $n$ scales from $10$ up to $1.34 \times 10^9$, doubling at each step. Item values $v$ and weights $w$ are randomly generated, with $W$ fixed at half of the total weight. Detailed results are reported in Table~\ref{tab:knapsack}.
\blue{Across the tested sizes, GFORS closely matches GRB in objective value while becoming faster as the instances grow. cuOpt is also competitive on small and medium instances, but its performance deteriorates at larger scales. For example, at $\mathrm{NNZ}=2.10\times10^7$ and $4.19\times10^7$, cuOpt returns objective values of only $4.40\times10^3$ and $4.41\times10^3$, respectively, compared with approximately $9.06\times10^7$ and $1.81\times10^8$ from GFORS and GRB. At $\mathrm{NNZ}=8.39\times10^7$, its incumbent drops further to $4.00\times10^1$, and it runs out of memory at all larger tested sizes. At $\mathrm{NNZ}=1.68\times10^8$, GFORS attains a better incumbent objective value than GRB, and at $\mathrm{NNZ}=3.36\times10^8$, GRB runs out of memory while GFORS still finds a high-quality solution in about $136$ seconds.} A distinguishing feature of the knapsack problem is its extremely high-dimensional decision space, which creates substantial memory demands for solution sampling. This cost can be mitigated by setting the batch size $k_b=1$ in Algorithm~\ref{algo:gfors}, enabling GFORS to handle much larger instances. The corresponding peak GPU memory usage is reported in the final column.

\subsubsection{Max-Cut: Unconstrained QP vs Linearized BIP}

Given a connected network $G=(V,E)$ with signed edge weights, the max-cut problem seeks a bipartition of $V$ whose induced cut maximizes the total weight of edges crossing the partition. This problem can be expressed directly as a quadratic program (QP) or reformulated as an integer program (IP) using the McCormick envelope:
\[
\begin{aligned}
  \text{(QP): } \quad
  & \max_{x \in \{0,1\}^V}
    \sum_{(i,j) \in E} w_{ij}\,(x_i + x_j - 2x_i x_j)
  &
  ~~
  \text{(IP): } \quad
  & \begin{aligned}[t]
      \max_{\substack{x \in \{0,1\}^V \\ y \in \{0,1\}^E}}
      & \sum_{(i,j) \in E} w_{ij}\,(x_i + x_j - 2y_{ij}) \\
      \text{s.t.} \quad
      & y_{ij} \leq x_i, \ \forall (i,j)\in E, \\
      & y_{ij} \leq x_j, \ \forall (i,j)\in E, \\
      & y_{ij} \geq x_i + x_j - 1, \ \forall (i,j)\in E.
    \end{aligned}
\end{aligned}
\]
% \[
%   \begin{aligned}
%     \text{(QP): } \quad & \max_{x \in \{0,1\}^V} \ \sum_{(i,j) \in E} w_{ij}\,(x_i + x_j - 2x_i x_j) 
%     & \quad 
%     \text{(IP): } \quad & \max_{\substack{x \in \{0,1\}^V \\ y \in \{0,1\}^E}} \ \sum_{(i,j) \in E} w_{ij}\,(x_i + x_j - 2y_{ij}) \\[0.5em]
%     & 
%     & \text{s.t. } & y_{ij} \leq x_i, \quad \forall (i,j)\in E, \\
%     & 
%     & & y_{ij} \leq x_j, \quad \forall (i,j)\in E, \\
%     & 
%     & & y_{ij} \geq x_i + x_j - 1, \quad \forall (i,j)\in E.
%   \end{aligned}
% \]
In our experiments, the number of vertices $n:=|V|$ ranges from $10$ to $20{,}480$. Graph topologies are generated randomly with density fixed at $0.5$, \blue{i.e., $\lfloor 0.5\binom{n}{2}\rfloor$ undirected edges are sampled uniformly without replacement, and disconnected instances are discarded.
Edge weights $w$ are sampled uniformly from the integers $\{-8,\ldots,10\}$. For the quadratic formulation, Gurobi is run with the setting \texttt{NonConvex=2}.} Results are reported in Table~\ref{tab:maxcut}.

\begin{table}[!tbp]
  \centering
  \blue{
 \scriptsize
\begin{tabular}{rrrrrrrrrr}
\toprule
Config & \multirow{2}{*}{NNZ}
& \multicolumn{4}{c}{Avg. Objective Value}
& \multicolumn{4}{c}{Avg. Runtime (sec.)} \\
\cmidrule(lr){3-6}\cmidrule(lr){7-10}
$n$ & & GRB & GFORS & cuOpt (IP) & GFORS (IP)
& GRB & GFORS & cuOpt (IP) & GFORS (IP) \\
\midrule
$10$ & $4.48e1$ & $4.12e1$ & $4.08e1$ & $4.12e1$ & $8.80e0$
& $0.01$ & $0.62$ & $0.77$ & $43.41$ \\
$20$ & $1.78e2$ & $1.70e2$ & $1.64e2$ & $1.70e2$ & $3.72e1$
& $0.10$ & $0.73$ & $0.73$ & $1471.68$ \\
$40$ & $7.41e2$ & $5.11e2$ & $4.10e2$ & $5.11e2$ & $2.10e2$
& $0.91$ & $0.95$ & $4.62$ & $159.24$ \\
$80$ & $3.01e3$ & $1.87e3$ & $1.59e3$ & $1.76e3$ & $3.23e2$
& $999.14$ & $0.84$ & $1805.32$ & $274.74$ \\
$160$ & $1.20e4$ & $6.07e3$ & $4.31e3$ & $5.43e3$ & $\mathrm{TL}$
& $1800.05$ & $1.18$ & $1801.62$ & $1800.01$ \\
$320$ & $4.82e4$ & $2.10e4$ & $1.41e4$ & $1.56e4$ & $2.90e2$
& $1800.07$ & $1.17$ & $1807.28$ & $1485.19$ \\
$640$ & $1.94e5$ & $7.46e4$ & $5.38e4$ & $5.76e4$ & $8.52e2$
& $1800.16$ & $1.18$ & $1800.59$ & $6.16$ \\
$1,280$ & $7.76e5$ & $2.72e5$ & $2.10e5$ & $1.17e5$ & $\mathrm{TL}$
& $1800.43$ & $1.17$ & $1809.31$ & $1800.01$ \\
$2,560$ & $3.10e6$ & $1.01e6$ & $8.32e5$ & $8.33e4$ & $4.94e3$
& $1800.95$ & $1.18$ & $1991.99$ & $17.13$ \\
$5,120$ & $1.24e7$ & $3.82e6$ & $3.30e6$ & $7.49e4$ & $8.66e3$
& $1802.87$ & $1.27$ & $2092.73$ & $71.71$ \\
$10,240$ & $4.97e7$ & $1.46e7$ & $1.31e7$ & $5.54e4$ & $2.06e4$
& $1842.38$ & $2.07$ & $107.22$ & $399.97$ \\
$20,480$ & $1.99e8$ & $5.66e7$ & $5.25e7$ & $\mathrm{OOM}$ & $\mathrm{OOM}$
& $1819.13$ & $5.15$ & $16.33$ & -- \\
$40,960$ & $7.95e8$ & $\mathrm{OOM}$ & $2.10e8$ & -- & --
& -- & $18.48$ & -- & -- \\
\bottomrule
\end{tabular}

  }
  \caption{Max-Cut: objective values and runtimes. \emph{GRB} and \emph{GFORS} solve the QP formulation, while \blue{\emph{cuOpt (IP)}} and \emph{GFORS (IP)} solve the IP formulation; \blue{$\mathrm{TL}$ and $\mathrm{OOM}$ indicate no feasible solution before the time or memory limit, respectively. Size $40{,}960$ is omitted for the IP formulation because its NNZ exceeds the $2\times10^9$ limit.}}
  \label{tab:maxcut}
\end{table}

% In our experiments, the number of vertices $n:=|V|$ ranges from $10$ to $20{,}480$. Graph topologies are generated randomly with density fixed at $0.5$, and edge weights $w$ are sampled uniformly from $[-8,10]$. Results are reported in Table~\ref{tab:maxcut}.

Both GRB and GFORS solve the QP formulation, while \blue{cuOpt (IP) and GFORS (IP) are applied to the IP reformulation. GFORS is within $4\%$ of GRB on the two smallest instances. The gap increases to about $33\%$ at intermediate sizes and then steadily narrows to about $7\%$ at the largest instance solved by both methods, while GFORS achieves shorter runtimes on medium and large instances. In comparison, cuOpt (IP) matches GRB on the three smallest instances but deteriorates rapidly as the problem size grows. For example, at $n=10{,}240$, it returns an objective value of only $5.54\times10^4$, compared with $1.31\times10^7$ for GFORS, and it runs out of memory at $n=20{,}480$. GFORS (IP) performs substantially worse than GFORS on the QP formulation and also fails on several instances. This stronger performance of GFORS on the QP formulation is consistent with the tighter relaxation provided by the QP model relative to the IP reformulation~\citep{fakhimi2023relaxations}. As suggested by Theorem~\ref{thm:optprob}, when a tighter relaxation yields a fractional solution closer to an optimal binary solution, the resulting smaller $\mu$ can strengthen the sampling guarantee and improve GFORS performance.}

\subsubsection{3D Assignment: Customized Sampling}
\label{sec:assign3d}

\begin{table}[!tbp]
  \centering
  \blue{
 \scriptsize
\begin{tabular}{lrrrrrrrrr}
\toprule
Config & \multirow{2}{*}{NNZ}
& \multicolumn{4}{c}{Avg. Objective Value}
& \multicolumn{4}{c}{Avg. Runtime (sec.)} \\
\cmidrule(lr){3-6}\cmidrule(lr){7-10}
$n$ & & GRB & cuOpt & Cust. & Def.
& GRB & cuOpt & Cust. & Def. \\
\midrule
$2$ & $1.60e1$ & $4.6$ & $4.6$ & $4.6$ & $4.6$
& $0.00$ & $0.31$ & $3.52$ & $1.11$ \\
$4$ & $1.60e2$ & $6.2$ & $6.2$ & $6.4$ & $8.0$
& $0.00$ & $0.33$ & $5.24$ & $1.52$ \\
$8$ & $1.41e3$ & $8.0$ & $8.0$ & $9.4$ & $16.4$
& $0.01$ & $0.41$ & $10.97$ & $6.44$ \\
$16$ & $1.18e4$ & $16.0$ & $16.0$ & $16.2$ & $\mathrm{TL}$
& $0.13$ & $2.88$ & $8.79$ & $1800.01$ \\
$32$ & $9.63e4$ & $32.0$ & $32.0$ & $32.0$ & $\mathrm{TL}$
& $0.64$ & $51.84$ & $22.39$ & $1800.01$ \\
$64$ & $7.78e5$ & $64.0$ & $64.0$ & $64.0$ & $\mathrm{TL}$
& $8.30$ & $216.75$ & $17.83$ & $1800.01$ \\
$128$ & $6.26e6$ & $128.0$ & $515.8$ & $128.0$ & $\mathrm{TL}$
& $75.75$ & $2072.74$ & $54.26$ & $1800.01$ \\
$256$ & $5.02e7$ & $256.0$ & $1333.8$ & $256.0$ & $\mathrm{TL}$
& $657.75$ & $2053.73$ & $250.16$ & $1800.04$ \\
$512$ & $4.02e8$ & $\mathrm{OOM}$ & $\mathrm{OOM}$ & $512.0$ & $\mathrm{TL}$
& $255.04$ & $14.18$ & $67.56$ & $1800.62$ \\
\bottomrule
\end{tabular}

  }
 \caption{3D Assignment: objective values and runtimes. \emph{Cust.} and \emph{Def.} denote GFORS with customized and default sampling, respectively; \blue{$\mathrm{TL}$ and $\mathrm{OOM}$ indicate no feasible solution before the time or memory limit, respectively.}}
  \label{tab:assign}
\end{table}

The 3D assignment problem generalizes the classic 2D assignment to a three-dimensional setting and is formulated as \eqref{eq:assign3d}. 
\blue{
We consider $n$ ranging from $2$ to $512$, with \(n^3\) binary variables and independent costs uniformly distributed on \(\{1,\ldots,10\}\). 
}
Even with the TU reformulation that eliminates the first two sets of constraints, $n$ equality constraints remain, limiting sampling efficiency. To address this, we implement the proposed GPU-friendly customized sampling procedure (Algorithm~\ref{alg:3dassign}), with results summarized in Table~\ref{tab:assign}.

Columns \emph{Cust.} and \emph{Def.} report results for customized and default sampling, respectively. The customized scheme yields a significant performance gain in both objective quality and runtime. \blue{With default sampling, only the three smallest instances yield feasible solutions, while all larger instances reach the time limit without one.} In contrast, with customized sampling, all instances are solved \blue{and GRB is matched on nearly all medium and large sizes, while requiring less time on the larger cases. cuOpt is competitive on small instances but degrades markedly at larger sizes and fails on the largest instance. At $n=512$, both GRB and cuOpt fail, whereas customized GFORS still finds the objective $512$ in under $70$ seconds.} These results demonstrate that customized sampling subroutines can significantly enhance the performance \blue{and scalability} of GFORS.

\paragraph{Effects of the First-Order Updates.}
\begin{table}[!bp]
  \centering
  \blue{
 \scriptsize
\setlength{\tabcolsep}{2.9pt}
\begin{tabular}{crrrrrrrrr}
\toprule
\multirow{2}{*}{Instance size}
& \multicolumn{3}{c}{10 LS steps}
& \multicolumn{3}{c}{50 LS steps}
& \multicolumn{3}{c}{100 LS steps} \\
\cmidrule(lr){2-4}\cmidrule(lr){5-7}\cmidrule(lr){8-10}
& Early mean & Late mean & Impr. (\%)
& Early mean & Late mean & Impr. (\%)
& Early mean & Late mean & Impr. (\%) \\
\midrule
64
& 123.00 & 99.30 & 19.27
& 78.97 & 73.80 & 6.54
& 70.93 & 68.30 & 3.71 \\

128
& 254.33 & 233.60 & 8.15
& 163.77 & 158.70 & 3.09
& 142.33 & 140.50 & 1.29 \\

256
& 527.07 & 502.47 & 4.67
& 329.47 & 324.87 & 1.40
& 287.53 & 284.93 & 0.90 \\

512
& 1073.53 & 1045.07 & 2.65
& 672.37 & 665.60 & 1.01
& 580.83 & 578.93 & 0.33 \\
\bottomrule
\end{tabular}

  }
  \caption{Effects of first-order distribution updates on the 3D assignment problem with customized sampling. For each number of local-search steps ($L$ in Algorithm~\ref{alg:3dassign}), the table reports the mean minimum sampled objective over iterations 200--250 (Early) and 450--500 (Late), together with the corresponding relative improvement (\%); positive values indicate improved (lower) objectives.}
  \label{tab:apquantile}
\end{table}

\blue{
A natural question is whether the improvement is driven by first-order guidance of the sampling distribution or primarily by the local search in Algorithm~\ref{alg:3dassign}. To isolate the effect of the first-order updates, we compare the minimum sampled objective values over two equal-length windows: iterations \(200\text{--}250\) (Early), before the update takes effect, and iterations \(450\text{--}500\) (Late), after meaningful updates have been made. For each instance size and local-search setting, the minimum objective is averaged within the two windows and across the five instances. Since the assignment problem is a minimization problem, we report
\(
100\times (v_{\mathrm{before}}-v_{\mathrm{after}})/v_{\mathrm{before}},
\)
so positive values indicate improvement after the first-order update.
Table~\ref{tab:apquantile} reports the results. The first-order updates improve the minimum sampled objective across every instance size and step setting, but the effect consistently diminishes as local search becomes more extensive. For example, the improvement ranges from \(2.65\%\) to \(19.27\%\) with \(10\) steps, compared with only \(0.33\%\) to \(3.71\%\) with \(100\) steps. The improvement also decreases as instance size grows, suggesting that larger instances require more optimization progress before the sampling distribution can be meaningfully refined. Overall, first-order updates provide greater guidance when local search is shallow, while stronger local search reduces its marginal contribution.}

\subsubsection{Facility Location: TU Reformulation}
\label{sec:fl}

This version of the facility location problem seeks the optimal set of facility openings and customer assignments that minimize total cost.
\[
  \begin{aligned}
    \min_{x \in \{0,1\}^{n_f},\, y \in \{0,1\}^{n_f \times n_c}} 
    &~ \iprod{f, x} + \iprod{c, y} \\
    \text{s.t. } 
    &~ \sum_{i \in F} y_{ij} = 1, \quad \forall j \in C, \\
    &~ y_{ij} \leq x_i, \quad \forall i \in F,~ j \in C.
  \end{aligned}
\]
The size parameters $(n_f,n_c)$ range from $(2,8)$ up to $(8192,32768)$, with
\blue{facility-opening costs sampled uniformly from the integers \(\{10,\ldots,21\}\) and customer-assignment costs from \(\{1,\ldots,11\}\).} 
\blue{Results are reported in Table~\ref{tab:fl}.}

\begin{table}[!tbp]
  \centering
  \blue{
 \scriptsize
\begin{tabular}{lrrrrrrrrrr}
\toprule
Config & \multirow{2}{*}{NNZ}
& \multicolumn{4}{c}{Avg. Objective Value}
& \multicolumn{4}{c}{Avg. Runtime (sec.)}
& \multirow{2}{*}{TU Time} \\
\cmidrule(lr){3-6}\cmidrule(lr){7-10}
$(n_f,n_c)$ & & GRB & cuOpt & TU & No TU
& GRB & cuOpt & TU & No TU & \\
\midrule
$(2,8)$ & $4.80e1$ & $61.4$ & $61.4$ & $61.4$ & $61.4$
& $0.00$ & $0.33$ & $0.75$ & $1.13$ & $0.00$ \\
$(4,16)$ & $1.92e2$ & $90.8$ & $90.8$ & $95.0$ & $92.0$
& $0.00$ & $0.33$ & $2.01$ & $2.96$ & $0.00$ \\
$(8,32)$ & $7.68e2$ & $133.6$ & $133.6$ & $140.2$ & $137.6$
& $0.01$ & $0.41$ & $8.05$ & $5.61$ & $0.00$ \\
$(16,64)$ & $3.07e3$ & $188.4$ & $188.4$ & $224.2$ & $244.8$
& $0.06$ & $1.69$ & $8.60$ & $6.91$ & $0.00$ \\
$(32,128)$ & $1.23e4$ & $296.0$ & $296.0$ & $545.2$ & $436.6$
& $1.72$ & $10.09$ & $7.05$ & $10.43$ & $0.00$ \\
$(64,256)$ & $4.92e4$ & $458.6$ & $1520.6$ & $1561.2$ & $952.2$
& $55.29$ & $23.81$ & $0.80$ & $11.30$ & $0.00$ \\
$(128,512)$ & $1.97e5$ & $767.0$ & $766.6$ & $3096.8$ & $\mathrm{TL}$
& $1800.25$ & $1801.31$ & $0.83$ & $1800.01$ & $0.01$ \\
$(256,1024)$ & $7.86e5$ & $1340.0$ & $1325.2$ & $6135.8$ & $\mathrm{TL}$
& $1800.51$ & $1802.67$ & $0.90$ & $1800.01$ & $0.04$ \\
$(512,2048)$ & $3.15e6$ & $2489.5$ & $2408.4$ & $12260.8$ & $\mathrm{TL}$
& $1800.91$ & $1815.27$ & $1.25$ & $1800.01$ & $0.16$ \\
$(1024,4096)$ & $1.26e7$ & $22498.8$ & $5735.0$ & $24586.2$ & $\mathrm{TL}$
& $1802.14$ & $1806.85$ & $2.37$ & $1800.02$ & $0.83$ \\
$(2048,8192)$ & $5.03e7$ & $80296.8$ & $49205.4$ & $49157.4$ & $\mathrm{TL}$
& $1186.78$ & $158.67$ & $9.54$ & $1800.09$ & $5.07$ \\
$(4096,16384)$ & $2.01e8$ & $160808.4$ & $\mathrm{OOM}$ & $98476.9$ & $\mathrm{TL}$
& $993.81$ & $50.97$ & $36.81$ & $1800.25$ & $29.53$ \\
$(8192,32768)$ & $8.05e8$ & $\mathrm{OOM}$ & $\mathrm{OOM}$ & $\mathrm{OOM}$ & $\mathrm{OOM}$
& -- & $22.88$ & -- & -- & -- \\
\bottomrule
\end{tabular}

  }
 \caption{Facility Location: objective values and runtimes. \emph{TU} and \emph{No TU} denote the implementations of GFORS \blue{with and without} TU reformulation; \emph{Avg. TU Time} denotes the mean wall-clock time spent on the TU reformulation; \blue{$\mathrm{TL}$ and $\mathrm{OOM}$ indicate no feasible solution before the time or memory limit, respectively.}}
  \label{tab:fl}
\end{table}

Columns \emph{TU} and \emph{No TU} correspond to GFORS implementations with and without the TU reformulation step (see Section~\ref{sec:tureform}). With a short execution time (see \emph{Avg. TU Time}), the TU reformulation significantly improves GFORS performance: \blue{all solved configurations are completed within about $37$ seconds, whereas the No-TU variant fails to find a feasible solution from $(128,512)$ onward.} 
In terms of solution quality, GRB \blue{and cuOpt} perform best on small- and medium-sized instances, while GFORS (TU) \blue{begins to outperform GRB at $(2048,8192)$ and remains better at $(4096,16384)$. At $(2048,8192)$, GFORS (TU) also slightly outperforms cuOpt. At $(4096,16384)$, cuOpt runs out of memory, while GFORS (TU) obtains a lower objective than GRB in $36.81$ seconds.} These results highlight that TU reformulation can greatly enhance GFORS's performance when it eliminates all equality constraints.

% Columns \emph{TU} and \emph{No TU} correspond to GFORS implementations with and without the TU reformulation step (see Section~\ref{sec:tureform}). With a short execution time (see \emph{Avg. TU Time}), the TU reformulation significantly improves GFORS performance: \blue{all solved configurations are completed within about $37$ seconds, whereas the No-TU variant fails to find a feasible solution from $(128,512)$ onward.} In terms of solution quality, GRB dominates on small- and medium-sized instances, but GFORS (TU) \blue{begins to outperform GRB at $(1024,4096)$. At this size, cuOpt obtains the best objective. However, GFORS (TU) achieves the best objective at $(2048,8192)$ and $(4096,16384)$ while requiring only $9.54$ and $36.81$ seconds, respectively.} These results highlight that TU reformulation can greatly enhance GFORS's performance when it eliminates all equality constraints.

% \subsubsection{TSP: Challenges from Weak Relaxation \& Sparse Feasible Set}
% \label{sec:tsp}

\subsection{ Incumbent Trajectories under Randomized Sampling}
In Algorithm~\ref{algo:gfors}, the parameter \(k_{\mathrm{int}}\) controls how frequently randomized sampling is invoked to update the incumbent, trading off per-iteration overhead against the breadth of sampling coverage.
\blue{
To assess iterative sampling, we examine incumbent trajectories across multiple problem classes, highlighting when improvements occur and how intermediate sampling contributes to solution quality.
}

\begin{figure}[!t]
  \centering
  \begin{subfigure}[t]{0.33\linewidth}
    \centering\includegraphics[width=\linewidth]{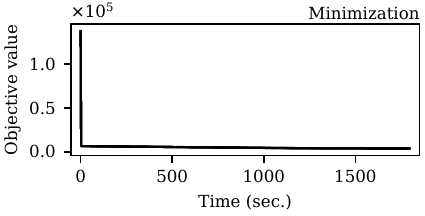}\caption{Set Cover $\sigma=0.99$}\label{fig:1}
  \end{subfigure}\hfill
  \begin{subfigure}[t]{0.33\linewidth}
    \centering\includegraphics[width=\linewidth]{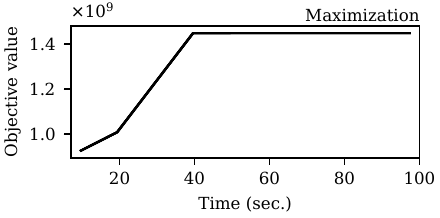}\caption{Knapsack}\label{fig:2}
  \end{subfigure}\hfill
  \begin{subfigure}[t]{0.33\linewidth}
    \centering\includegraphics[width=\linewidth]{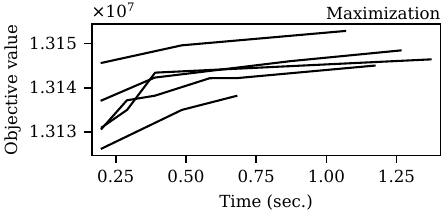}\caption{Max Cut QP}\label{fig:3}
  \end{subfigure}

  \vspace{0pt}

  \begin{subfigure}[t]{0.34\linewidth}
    \centering\includegraphics[width=\linewidth]{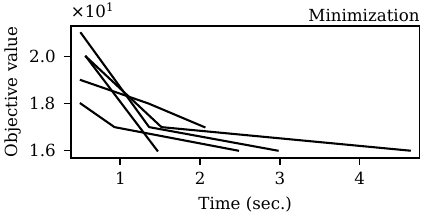}\caption{3D Assignment Cust.}\label{fig:4}
  \end{subfigure}\quad
  \begin{subfigure}[t]{0.34\linewidth}
    \centering\includegraphics[width=\linewidth]{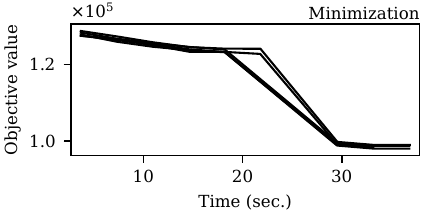}\caption{Facility Location TU}\label{fig:5}
  \end{subfigure}

  \caption{Iterative incumbent objective value versus runtime (seconds) by problem class \blue{without aggregation}. Across most classes, randomized sampling yields rapid early improvements followed by slower refinement, \blue{highlighting its ability to identify strong incumbents early in the search.}}

  \label{fig:traj}
\end{figure}

Specifically, for each class, we select five randomly generated instances and plot the incumbent objective value versus runtime in seconds \blue{without aggregation}. We focus on relatively large instance settings to ensure that each run performs multiple sampling subroutines, making the effect of iterative sampling visible. The resulting trajectories are shown in Figure~\ref{fig:traj}.

For set cover and knapsack, the selected instances show a similar pattern: high-quality solutions are obtained in early iterations, while refinements become slow in the late stage. For facility location, however, the most significant improvements occur relatively later in the execution, followed by a refinement phase. These trends are less apparent for max-cut and 3D assignment, potentially due to their shorter execution time. Overall, randomized sampling supports a refinement scheme that finds good incumbents quickly and then improves them gradually.

\subsection{\blue{QOBLIB Maximum Independent Set (MIS) Instances}}
\blue{
This subsection evaluates GFORS on public instances from the QOBLIB repository~\citep{koch2026quantumoptimizationbenchmarkinglibrary,koch2025quantumoptimizationbenchmarkinglibrary}. We focus on the MIS problem, which provides a natural binary integer programming benchmark with well-established hard instances and a standard linear formulation. 
Among the 50 MIS instances in QOBLIB, we select the ten largest by NNZ, providing a challenging test set that spans the largest instances in the repository.
We consider the following two formulations:
\[
\begin{aligned}
  \text{(QP): } \quad
  & \max_{x \in \{0,1\}^V}
  \left\{
    \sum_{i\in V} x_i
    -
    M \sum_{(i,j)\in E} x_i x_j
  \right\}
  &
  ~~
  \text{(IP): } \quad
  & \begin{aligned}[t]
      \max_{x \in \{0,1\}^V} \quad
      & \sum_{i\in V} x_i \\
      \text{s.t.} \quad
      & x_i + x_j \leq 1,
      \ \forall (i,j)\in E.
    \end{aligned}
\end{aligned}
\]
We set \(M:=2\), since any \(M>1\) guarantees that the global optimizers of the penalized QUBO coincide with those of the original MIS problem~\citep{lucas2014ising}. For each instance, Gurobi and cuOpt (both with presolve enabled) solve the IP formulation, while GFORS uses the QP formulation. The latter is better suited to GFORS because, on dense graphs, the IP relaxation often has an optimum near \(x=0.5\). With GFORS initialized at the same point, the resulting first-order trajectory can provide little guidance for sampling. We also use the following customized sampling and terminal refinement along with QP-based GFORS.
}

\blue{
\paragraph{Customized Sampling for MIS.}
At each sampling checkpoint, given the current fractional solution \(x\), we assign vertex \(i\) the randomized priority
\(
p_i=\log\!\bigl(\max\{x_i,\epsilon\}\bigr)+ g_i,
\)
where \(\epsilon>0\) prevents numerical issues and \(g_i\) is an independent Gumbel variate. Vertices are processed in decreasing priority: the highest-priority active vertex is selected and its neighbors are deactivated, yielding a maximal independent set. Sampling is performed in parallel on the GPU, and up to \(500\) best unique candidates are retained throughout the GFORS run.
}

\begin{table}[!tbp]
  \centering
  \blue{
 \centering
\scriptsize
\setlength{\tabcolsep}{3pt}
\resizebox{\linewidth}{!}{%
\begin{tabular}{llrrrrrrrrr}
\toprule
\multirow{2}{*}{Instance}
  & \multirow{2}{*}{$(n,m)$}
  & \multirow{2}{*}{IP NNZ}
  & \multicolumn{4}{c}{Objective Value}
  & \multicolumn{4}{c}{Time (sec.)} \\
\cmidrule(lr){4-7}\cmidrule(lr){8-11}
  & & & GRB & cuOpt & GFORS & GFORS+TR
  & GRB & cuOpt & GFORS & GFORS+TR \\
\midrule
frb53-24-1 & $(1272{,}\ 94227)$ & $1.88e5$ & 50 & 51 & 40 & 46 & 31.46 & 1.79 & 0.77 & +0.15 \\
socfb-trinity100 & $(2613{,}\ 111996)$ & $2.24e5$ & 499 & 499 & 474 & 487 & 69.57 & 0.45 & 1.47 & +0.42 \\
brock800-1 & $(800{,}\ 112095)$ & $2.24e5$ & 19 & 21 & 18 & 19 & 43.89 & 0.28 & 0.74 & +0.08 \\
frb59-26-2 & $(1534{,}\ 126163)$ & $2.52e5$ & 55 & 56 & 45 & 51 & 63.96 & 0.28 & 0.79 & +0.20 \\
p\_hat1500-3 & $(1500{,}\ 277006)$ & $5.54e5$ & 93 & 94 & 88 & 92 & 192.36 & 0.39 & 0.83 & +1.67 \\
sorrell4 & $(2048{,}\ 504451)$ & $1.01e6$ & 23 & 24 & 19 & 21 & 0.02 & 0.70 & 0.76 & +0.40 \\
frb100-40 & $(4000{,}\ 572774)$ & $1.15e6$ & 76 & 92 & 74 & 84 & 1282.75 & 0.61 & 0.84 & +1.37 \\
p\_hat1500-1 & $(1500{,}\ 839327)$ & $1.68e6$ & 11 & 11 & 10 & 11 & 552.26 & 0.72 & 0.76 & +0.50 \\
keller6 & $(3361{,}\ 1026582)$ & $2.05e6$ & 44 & 57 & 44 & 51 & 1777.25 & 0.89 & 0.80 & +2.22 \\
C4000-5 & $(4000{,}\ 3997732)$ & $8.00e6$ & 12 & 17 & 15 & 15 & $\mathrm{TL}$ & 2.98 & 0.81 & +2.73 \\
\bottomrule
\end{tabular}}

  }
  \caption{\blue{MIS instances from QOBLIB. GRB and cuOpt are run for the full 1,800-second time limit without closing the optimality gap. Their Time columns report the first time their objective reaches or exceeds the GFORS objective. The GFORS Time column reports its solve time, while GFORS+TR reports the additional GPU time for terminal refinement. $\mathrm{TL}$ indicates that the GFORS objective was not reached within 1,800 seconds.}}
  \label{tab:mis}
\end{table}

\blue{
\paragraph{Terminal Refinement for MIS.}
After GFORS terminates, each retained maximal independent set is refined on the GPU in parallel using batched \(1\)-to-\(2\) exchanges~\citep{andrade2012fast}. Each accepted exchange removes one selected vertex and inserts two nonadjacent unselected vertices for which the removed vertex is their only selected neighbor, thereby increasing the cardinality by one. Newly admissible vertices are then greedily added to restore maximality. The procedure terminates when no improving exchange is found or the step limit of refinement is reached. All candidates are then verified, and the largest valid independent set is returned.
}

\blue{
The results are reported in Table~\ref{tab:mis}. Both GRB and cuOpt run for the full 1,800-second limit without closing the optimality gap. Their Time columns report the first time at which their incumbent reaches or exceeds the GFORS objective before terminal refinement. The GFORS Time column reports its solve time, while GFORS+TR reports the additional GPU time required by terminal refinement.
Overall, cuOpt attains the best objective on most instances and is the fastest to reach the GFORS objective in seven of the ten cases. GFORS nevertheless reaches its incumbent in about one second across all instances and is faster than cuOpt on \texttt{frb53-24-1} and on the two largest instances, \texttt{keller6} and \texttt{C4000-5}. Gurobi is generally much slower to reach the same objective and fails to do so within the time limit on \texttt{C4000-5}. 
Terminal refinement improves the GFORS solutions on nearly all instances at modest additional GPU cost, narrowing the remaining objective gap and matching or surpassing Gurobi on several large instances.
}

\section{Conclusion}
\label{sec:conc}
Across the tested problem classes, our study shows that a GPU-native combination of a PDHG-style first-order routine with repeated batched sampling delivers strong time-to-incumbent on large BIPs. The framework runs end-to-end on GPUs with minimal synchronization and leverages TU reformulations and customized sampling to enhance sampling efficiency. 
% It provides near-stationarity guarantees for the first-order component together with probabilistic bounds on sampled solutions. 
\blue{It provides an ${\cal O}(N^{-1/2})$ residual convergence guarantee for the first-order component in the convex setting and a finite-horizon residual estimate in the general nonconvex setting, together with probabilistic bounds on sampled solutions.}
In practice, \blue{Gurobi and cuOpt} remain stronger on small to medium instances, while GFORS offers complementary scalability and increasingly favorable performance on large instances under tight time limits.

% In practice, exact solvers such as Gurobi remain stronger on small to medium instances, while GFORS offers complementary scalability on large instances under tight time limits.

We highlight several directions for future work. First, integrating lower-bounding mechanisms is a priority to strengthen optimality guarantees. Second, the current design assumes explicitly enumerated constraints, making it ill-suited to formulations that rely on dynamic separation or an exponential number of cuts. Third, extending beyond binary variables to general integer and mixed-integer programs requires further development. 
Fourth, while the current implementation emphasizes an end-to-end, GPU-native workflow, additional systems-level optimizations may further improve performance.
Finally, from an architectural perspective, multi-GPU and distributed variants are a natural next step to expand problem scale and improve time-to-incumbent performance.

\newpage
\bibliographystyle{plainnat}
% \nocite{*}
\bibliography{bibi/refs}

\newpage
\appendix
\section{Mathematical Proofs}
% \section{Mathematical Proofs}
\label{sec:proofs}
\propconv*

\pfstart
(a) The optimality condition of \eqref{eq:pdhg-y} is
\[
0 \in \frac{y_k-y_{k-1}}{\tau_2}-\nabla_y \hL(\bar x_{k-1}, y_k) + \partial h_2(y_k),
\]
which together with \eqref{eq:dy} and the definition of $y^k$ in \eqref{def:yk} implies that
    \begin{align*}
        0 &\in \frac{y^k-y^{k-1}}{\tau_2} - K(\bar x_k - x_k - \bar x_{k-1} +  x_{k-1})- (K \bar x_{k-1} +r) + \partial h_2(y_k) \\
        &= \frac{y^k-y^{k-1}}{\tau_2} - (K x_k + r) + \partial h_2(y_k) 
        \stackrel{\eqref{eq:dy}}= \frac{y^k-y^{k-1}}{\tau_2} - \nabla_y \hL(x_k, y_k) + \partial h_2(y_k). 
    \end{align*}
    Hence, the first inclusion in \eqref{incl:IPP} follows.
    The optimality condition of \eqref{eq:pdhg-x} is
\[
 \frac{x_{k-1} - x_k}{\tau_1} \in \partial [\ell_{\hL(\cdot,y_k)}(\cdot;x_{k-1}) +h_1(\cdot)](x_k),
\]
which together with the convexity of $\hL(\cdot,y_k)$ further implies that for every $x \in \dom h_1$,
\begin{align}
    &\hL(x,y_k)+h_1(x)\ge \ell_{\hL(\cdot,y_k)}(x;x_{k-1}) + h_1(x) \nn \\
    \ge& \ell_{\hL(\cdot,y_k)}(x_k;x_{k-1}) + h_1(x_k) + \frac{1}{\tau_1}\iprod{x_{k-1} - x_k, x-x_k} \nn \\
    =& \hL(x_k,y_k) - \varepsilon_k + h_1(x_k) + \frac{1}{\tau_1}\iprod{x_{k-1} - x_k, x-x_k}.\label{ineq:x0}
\end{align}
Hence, the second inclusion in \eqref{incl:IPP} follows.

(b) It follows from the definition of $\varepsilon_k$ in \eqref{def:epsk} and the assumption that $\hL(\cdot,y_k)$ is $L$-smooth that
    \[
    \varepsilon_k \le \frac L2 \|x_k-x_{k-1}\|^2  \le \frac{1}{8 \tau_1} \|x_k-x_{k-1}\|^2,
    \]
    where the second inequality is due to the assumption that $\tau_1 \le 1/(4L)$. Hence, the first inequality in \eqref{ineq:ek} holds. Plugging this inequality into \eqref{ineq:x0}, we have
    \begin{equation}\label{ineq:xk}
        \hL(x,y_k) + h_1(x) - \hL(x_k,y_k) - h_1(x_k)  \ge \frac{3}{8\tau_1}\|x_k-x_{k-1}\|^2 + \frac{1}{2\tau_1}\|x-x_k\|^2 - \frac{1}{2\tau_1}\|x-x_{k-1}\|^2.
    \end{equation}
    Using \eqref{eq:pdhg-bx} and \eqref{def:yk}, we obtain
	\[
	   \frac{1}{\tau_2}\|y^k-y_k\|^2 
        \stackrel{\eqref{def:yk}}= \tau_2\|K(\bar{x}_k- x_k)\|^2 
        \stackrel{\eqref{eq:pdhg-bx}}= \tau_2\|K(x_k- x_{k-1})\|^2 
        \le \tau_2 \|K\|^2 \|x_k-x_{k-1}\|^2. 
	\]
    Hence, the second inequality in \eqref{ineq:ek} follows from the assumption that $\tau_1 \tau_2 \|K\|^2 \le 1/4$.

(c) The first inclusion in \eqref{incl:IPP} implies that for every $y \in \dom h_2$,
    \begin{align*}
        &-\hL(x_k,y)+h_2(y) \ge -\hL(x_k,y_k) + h_2(y_k) +\frac{1}{\tau_2}\iprod{y^{k-1}-y^k, y-y_k} \\
        =& -\hL(x_k,y_k) + h_2(y_k) + \frac{1}{2\tau_2}\|y-y^k\|^2 -\frac{1}{2\tau_2}\|y-y^{k-1}\|^2 + \frac{1}{2\tau_2}\|y_k-y^{k-1}\|^2 - \frac{1}{2\tau_2}\|y_k-y^k\|^2.
    \end{align*}
    Plugging the second inequality in \eqref{ineq:ek} into the above inequality yields
    \begin{align}
        &-\hL(x_k,y)+h_2(y) +\hL(x_k,y_k) - h_2(y_k) \nn \\
    \ge&\frac{1}{2\tau_2}\left(\|y-y^k\|^2 -\|y-y^{k-1}\|^2 + \|y_k-y^{k-1}\|^2\right) - \frac{1}{8\tau_1} \|x_k-x_{k-1}\|^2. \label{ineq:y0}
    \end{align}
    Summing \eqref{ineq:xk} and \eqref{ineq:y0} gives
    \begin{align*}
        \Phi(x,y_k) - \Phi(x_k,y)
        &\ge \frac{1}{4\tau_1}\|x_k-x_{k-1}\|^2 + \frac{1}{2\tau_1}\|x-x_k\|^2 -\frac{1}{2\tau_1}\|x-x_{k-1}\|^2 \\
        & + \frac{1}{2\tau_2}\|y-y^k\|^2 -\frac{1}{2\tau_2}\|y-y^{k-1}\|^2 +\frac{1}{2\tau_2}\|y_k-y^{k-1}\|^2.
    \end{align*}
    Taking $(x,y)=(x_*,y_*)$ as the saddle point of $\Phi(\cdot,\cdot)$ and using the fact that $\Phi(x_*,y_k) \le \Phi(x_k,y_*)$, we obtain
    \[
    \frac{1}{4\tau_1} \|x_k-x_{k-1}\|^2+\frac{1}{2\tau_2}\|y_k-y^{k-1}\|^2 \le \frac{1}{2\tau_1}\|x_*-x_{k-1}\|^2 - \frac{1}{2\tau_1}\|x_*-x_k\|^2 + \frac{1}{2\tau_2}\|y_*-y^{k-1}\|^2 -  \frac{1}{2\tau_2}\|y_*-y^k\|^2.
    \]
    Finally, \eqref{ineq:IPP} follows by summing the above inequality over iterations.
\pfend

\thmconv*
\pfstart
Let $k_0$ be the index such that
\[
k_0 =\argmin\left\{ \frac{1}{4\tau_1} \|x_k-x_{k-1}\|^2+\frac{1}{2\tau_2} \|y_k-y^{k-1}\|^2: k=1, \ldots, N\right\}.
\]
It follows from Proposition \ref{prop:convex} (a) that \eqref{con:incl} holds.
Using Proposition \ref{prop:convex} (c), we have
	\begin{equation}\label{ineq:Delta}
	    \frac{1}{4\tau_1} \|x_{k_0}-x_{k_0-1}\|^2+\frac{1}{2\tau_2} \|y_{k_0}-y^{k_0-1}\|^2 \stackrel{\eqref{ineq:IPP}}\le \frac{\|x_*-x_0\|^2}{2\tau_1 N} + \frac{\|y_*-y_0\|^2}{2\tau_2 N} = \frac{\Delta_0}{N}.
	\end{equation}
    Hence, it follows from the first inequality in Proposition \ref{prop:convex} (b) that
    \[
    \|r^x_{k_0}\|\stackrel{\eqref{incl:IPP}}=\frac{\|x_{k_0}-x_{k_0-1}\|}{\tau_1} \le \frac{2\sqrt{\Delta_0}}{\sqrt{\tau_1 N}}, \quad \varepsilon_{k_0} \stackrel{\eqref{ineq:ek}}\le \frac{1}{8 \tau_1} \|x_{k_0}-x_{k_0-1}\|^2 \le \frac{\Delta_0}{2N}.
    \]
    We have thus proved the last two inequalities in \eqref{con:ineq}.
    Using the Cauchy-Schwarz inequality and the second inequality in \eqref{ineq:ek}, we have
    \[
    \frac{1}{\tau_2}\|y^k-y^{k-1}\|^2 \le 3 \left(\frac{1}{\tau_2}\|y^k-y_k\|^2 + \frac{1}{2\tau_2}\|y_k-y^{k-1}\|^2\right) \stackrel{\eqref{ineq:ek}}\le 3 \left(\frac{1}{4\tau_1}\|x_k-x_{k-1}\|^2 + \frac{1}{2\tau_2}\|y_k-y^{k-1}\|^2\right),
    \]
    which together with \eqref{ineq:Delta} implies that
    \[
    \|r^y_{k_0}\| \stackrel{\eqref{incl:IPP}}= \frac{\|y^{k_0}-y^{k_0-1}\|}{\tau_2} \le \frac{\sqrt{3\Delta_0}}{\sqrt{\tau_2 N}}.
    \]
    Therefore, we have thus proved the first inequality in \eqref{con:ineq}.
\pfend

\propsum*
\pfstart
    It follows from \eqref{eq:pdhg-x} that for every $x$,
    \begin{align*}
	& \ell_{\hL(\cdot,y_k)}(x;x_{k-1}) + h_1(x) + \frac{1}{2\tau_1}\|x-x_{k-1}\|^2 - \frac{1}{2\tau_1}\|x-x_k\|^2 \\
	\ge& \ell_{\hL(\cdot,y_k)}(x_k;x_{k-1}) + h_1(x_k) + \frac{1}{2\tau_1}\|x_k-x_{k-1}\|^2 \\
	\ge& \hL(x_k,y_k) + h_1(x_k) - \frac{L}{2}\|x_k-x_{k-1}\|^2 + \frac{1}{2\tau_1}\|x_k-x_{k-1}\|^2.
    \end{align*}
    Taking $x=x_{k-1}$ yields
\[
    \hL(x_{k-1},y_k)+ h_1(x_{k-1}) \ge
\hL(x_k,y_k) + h_1(x_k) + \left(\frac{1}{\tau_1} - \frac{L}{2}\right)\|x_k-x_{k-1}\|^2.
\]
Noting $h_2$ is an indicator function, and thus $h_2(y_{k-1})=h_2(y_k)=0$.
This observation and the above inequality thus imply that
\begin{align*}
    \Phi(x_{k-1},y_{k-1}) &- \Phi(x_k,y_k) = \hL(x_{k-1},y_{k-1}) - \hL(x_k,y_k) + h_1(x_{k-1}) - h_1(x_k) \\
    &= \hL(x_{k-1},y_{k-1}) - \hL(x_{k-1},y_k) + \hL(x_{k-1},y_k) - \hL(x_k,y_k) + h_1(x_{k-1}) - h_1(x_k) \\
    &\stackrel{\eqref{eq:dy}}\ge \iprod{K x_{k-1} + r, y_{k-1} - y_k} + \left(\frac{1}{\tau_1} - \frac{L}{2}\right)\|x_k-x_{k-1}\|^2.
\end{align*}
Rearranging the terms gives
\begin{equation}\label{ineq:x}
    \left(\frac{1}{\tau_1} - \frac{L}{2}\right)\|x_k-x_{k-1}\|^2 \le \Phi(x_{k-1},y_{k-1}) - \Phi(x_k,y_k) + \iprod{K x_{k-1} + r, y_k - y_{k-1}}.
\end{equation}
Similarly, it follows from \eqref{eq:pdhg-y} that for every $y$,
    \begin{align*}
		 -\hL(\bar x_{k-1},y)+h_2(y)+\frac{1}{2 \tau_2}\|y-y_{k-1}\|^2-\frac{1}{2 \tau_2}\|y-y_k\|^2 
		 \ge -\hL(\bar x_{k-1},y_k) + h_2(y_k) +\frac{1}{2 \tau_2}\|y_k-y_{k-1}\|^2.
	\end{align*}
Taking $y=y_{k-1}$ yields
\[
-\hL(\bar x_{k-1},y_{k-1})+h_2(y_{k-1}) 
		 \ge -\hL(\bar x_{k-1},y_k) + h_2(y_k) +\frac{1}{\tau_2}\|y_k-y_{k-1}\|^2.
\]
Noting $h_2(y_{k-1})=h_2(y_k)=0$ and using \eqref{eq:dy}, the above inequality becomes
\[
    \iprod{K \bar x_{k-1} +r, y_k -y_{k-1}} \stackrel{\eqref{eq:dy}}= \hL(\bar x_{k-1},y_k)-\hL(\bar x_{k-1},y_{k-1})
		 \ge \frac{1}{\tau_2}\|y_k-y_{k-1}\|^2.
\]
Summing the above inequality and \eqref{ineq:x}, we have
\begin{align*}
    &\left(\frac{1}{\tau_1} - \frac{L}{2}\right)\|x_k-x_{k-1}\|^2 + \frac{1}{\tau_2}\|y_k-y_{k-1}\|^2 - [\Phi(x_{k-1},y_{k-1}) - \Phi(x_k,y_k)] \\ 
    \le& 2\iprod{K x_{k-1} + r, y_k - y_{k-1}} + \iprod{K (\bar x_{k-1} - x_{k-1}), y_k -y_{k-1}} \\
    \stackrel{\eqref{eq:pdhg-bx}}=& 2\iprod{K x_{k-1} + r, y_k - y_{k-1}} + \iprod{K (x_{k-1} - x_{k-2}), y_k -y_{k-1}} \\
    \le& 2 (\|K\| \|x_{k-1}\| + \|r\|) \|y_k - y_{k-1}\| + \|K\| \|x_{k-1} - x_{k-2}\| \|y_k -y_{k-1}\|,
\end{align*}
where the last inequality is due to the Cauchy-Schwarz inequality.
It thus follows from Young's and \eqref{ineq:tau} that
\begin{align*}
    &\left(\frac{1}{\tau_1} - \frac{L}{2}\right)\|x_k-x_{k-1}\|^2 + \frac{1}{\tau_2}\|y_k-y_{k-1}\|^2 - [\Phi(x_{k-1},y_{k-1}) - \Phi(x_k,y_k)] \\ 
    \le& 4 \tau_2 \|K\|^2 \|x_{k-1}\|^2 + 4\tau_2 \|r\|^2 + \frac{1}{2\tau_2}\|y_k - y_{k-1}\|^2 + \tau_2 \|K\|^2 \|x_{k-1} - x_{k-2}\|^2 + \frac{1}{4\tau_2}\|y_k -y_{k-1}\|^2 \\
    \stackrel{\eqref{ineq:tau}}\le& \frac2{\tau_1} \|x_{k-1}\|^2 + 4\tau_2 \|r\|^2 + \frac{3}{4\tau_2}\|y_k - y_{k-1}\|^2 + \frac1{2\tau_1} \|x_{k-1} - x_{k-2}\|^2.
\end{align*}
Rearranging the terms and summing the resulting inequality over iterations, we have
\[
\left(\frac{1}{2\tau_1} - \frac L2\right)\sum_{k=1}^N \|x_k-x_{k-1}\|^2 + \frac{1}{4\tau_2}\sum_{k=1}^N \|y_k-y_{k-1}\|^2 \le \Delta\Phi_N + \frac2{\tau_1} \sum_{k=1}^N\|x_{k-1}\|^2 + 4\tau_2 \|r\|^2 N,
\]
where we use the fact that $x_0=x_{-1}$.
\pfend

\thmweakconv*
\pfstart
It follows from \eqref{eq:pdhg-x} that
    \[
0\in \nabla_x \hL(x_{k-1},y_k) + \partial h_1(x_k) + \frac{x_k -x_{k-1}}{\tau_1},
\]
which together with the definition of $s_k^x$ implies that
\[
s_k^x = \frac{x_{k-1} - x_k}{\tau_1} + \nabla_x \hL(x_k,y_k) - \nabla_x \hL(x_{k-1},y_k)\in \nabla_x \hL(x_k,y_k) + \partial h_1(x_k).
\]
Similarly, it follows from \eqref{eq:pdhg-y} that
\[
0\in - \nabla_y \hL(\bar x_{k-1},y_k) + \partial h_2(y_k) + \frac{y_k -y_{k-1}}{\tau_2},
\]
which together with the definition of $s_k^y$ implies that
\[
s_k^y = \frac{y_{k-1} -y_k}{\tau_2} - \nabla_y \hL(x_k,y_k) + \nabla_y \hL(\bar x_{k-1},y_k) \in - \nabla_y \hL(x_k,y_k) + \partial h_2(y_k).
\]
Thus, we conclude that \eqref{con:incl1} holds.
Moreover, the definition of $s_k^x$, \eqref{ineq:L}, and \eqref{ineq:tau} yield that 
\begin{equation}\label{ineq:skx}
    \|s_k^x\| \le \left(\frac{1}{\tau_1} + L \right) \|x_{k-1} - x_k\| \stackrel{\eqref{ineq:tau}}\le \frac{3}{2\tau_1} \|x_{k-1} - x_k\|.
\end{equation}
The definition of $s_k^y$, \eqref{eq:dy}, and \eqref{eq:pdhg-bx} imply that 
\begin{align*}
    \|s_k^y\| &\le \frac{1}{\tau_2} \|y_k -y_{k-1}\| + \|K\| \|x_k - \bar x_{k-1}\| \\
    &\stackrel{\eqref{eq:pdhg-bx}}\le \frac{1}{\tau_2} \|y_k -y_{k-1}\| + \|K\| (\|x_k - x_{k-1}\| + \|x_{k-1} - x_{k-2}\|).
\end{align*}
It follows from the fact that $x_0=x_{-1}$ that
\begin{align*}
    \sum_{k=1}^N \|s_k^y\| &\le \frac{1}{\tau_2} \sum_{k=1}^N \|y_k -y_{k-1}\| + \|K\| \sum_{k=1}^N (\|x_k - x_{k-1}\| + \|x_{k-1} - x_{k-2}\|) \\
    &\le \frac{1}{\tau_2} \sum_{k=1}^N \|y_k -y_{k-1}\| + 2\|K\| \sum_{k=1}^N\|x_k - x_{k-1}\|
\end{align*}
Combining this inequality with \eqref{ineq:skx}, we obtain
\begin{equation}\label{ineq:sxy}
    \sum_{k=1}^N (\|s_k^x\| + \|s_k^y\|) \le \frac{1}{\tau_2} \sum_{k=1}^N \|y_k -y_{k-1}\| + \left( 2\|K\| + \frac{3}{2\tau_1} \right)\sum_{k=1}^N\|x_k - x_{k-1}\|.
\end{equation}
Using the Cauchy-Schwarz inequality, we have
\begin{align*}
    &\left( \frac{1}{\tau_2} \sum_{k=1}^N \|y_k -y_{k-1}\| + \left( 2\|K\| + \frac{3}{2\tau_1} \right)\sum_{k=1}^N\|x_k - x_{k-1}\| \right)^2 \\
    \le & \left[\left( 2\|K\| + \frac{3}{2\tau_1} \right)^2 \tau_1 N + \frac{N}{\tau_2} \right] \left[ \frac{1}{\tau_1} \sum_{k=1}^N \|x_k-x_{k-1}\|^2 + \frac{1}{\tau_2}\sum_{k=1}^N \|y_k-y_{k-1}\|^2 \right].
\end{align*}
It thus follows from Proposition \ref{prop:sum} that
\begin{align*}
    &\frac{1}{\tau_2} \sum_{k=1}^N \|y_k -y_{k-1}\| + \left( 2\|K\| + \frac{3}{2\tau_1} \right)\sum_{k=1}^N\|x_k - x_{k-1}\| \\
    \stackrel{\eqref{ineq:sum}}\le & 2\sqrt{N} \left[\left( 2\|K\| + \frac{3}{2\tau_1} \right)^2 \tau_1 + \frac{1}{\tau_2} \right]^{1/2} \left[\Delta \Phi_N + \frac2{\tau_1} \sum_{k=1}^N\|x_{k-1}\|^2 + 4\tau_2 \|r\|^2 N \right]^{1/2}.
\end{align*}
Finally, \eqref{bound:xy} immediately follows by plugging the above inequality into \eqref{ineq:sxy}.
\pfend

\smpinter*
\pfstart
  The following identity holds for every $p \in [0,1]^n$ 
  $$\sum_{\hat x \in \{0,1\}^n} \pi(p, \hat x) = 1,$$
since $\pi(p, \cdot)$ is a probability distribution over $\{0,1\}^n$.
Without loss of generality, we focus on expanding the last entry of $x_p:=\sum_{\hat x \in \{0,1\}^n} \pi(p, \hat x)\hat x$ as follows, where $\bar x \in \{0,1\}^{n-1}$ and $\bar p \in [0,1]^{n-1}$ denote the projections of $\hat x$ and $p$ onto the first $n-1$ entries:
  \begin{align*}
    (x_p)_n &= \sum_{\{(\bar x, 1) \mid \bar x \in \{0,1\}^{n-1}\}} \pi(p, (\bar x, 1)) \cdot 1 + \sum_{\{(\bar x, 0) \mid \bar x \in \{0,1\}^{n-1}\}} \pi(p, (\bar x, 0)) \cdot 0\\
    &= \sum_{\bar x \in \{0,1\}^{n-1}} \left( \prod_{i \in [n-1]} p_i^{\bar x_i}(1-p_i)^{1-\bar x_i} \right) \, p_n \\[6pt]
&= p_n \sum_{\bar x \in \{0,1\}^{n-1}} \prod_{i \in [n-1]} p_i^{\bar x_i}(1-p_i)^{1-\bar x_i}\\
            % &= \sum_{\bar x \in \{0,1\}^{n-1}} \left(\prod_{i \in [n-1]}((1-p_i) + \bar x_i(2p_i - 1))\right)((1-p_n)+ 1\cdot (2p_n - 1))\\
            % &= p_n\sum_{\bar x \in \{0,1\}^{n-1}} \left(\prod_{i \in [n-1]}((1-p_i) + \bar x_i(2p_i - 1))\right)\\
            &= p_n \sum_{\bar x \in \{0,1\}^{n-1}} \pi(\bar p, \bar x)\\
            &=p_n.
  \end{align*}
  The last equality is due to the first identity applied to the case of $n-1$.
\pfend

\optprob*
\pfstart
    For any random vector $x \in \{0,1\}^n$ \blue{is distributed according to the product Bernoulli distribution parameterized by} $p$, $\|x - x^*\|_1 \leq \delta/L_f$ ensures $f(x) \leq f(x^*) + \delta$. Thus, it suffices to show
$$\Pr(\|x - x^*\|_1 > \delta/L_f) \leq \exp\left(\tfrac{\delta}{L_f}\bigl(1 - \log\tfrac{\delta}{\mu L_f}\bigr) - \mu \right).$$
Since $x$ is entrywise Bernoulli, the distance $D:=\|x-x^*\|_1 = \sum_{i \in [n]}\mathbf 1(x_i \neq x^*_i)$ follows the Poisson binomial distribution, which adopts the following Chernoff right-tail bound when $\mu \leq t$
$$\Pr(D > t) \leq \exp\left(t(1-\log \frac{t}{\mu}) - \mu \right).$$
This proves the bound of $\Pr(\|x-x^*\|_1 > \delta / L_f)$, and the claim follows from the sample size $k$.
\pfend

\smpfeas*
\pfstart
  Let $\gamma^j(p):= \iprod{A_j, p} - b_j$ and $\gamma_+^j(p):= \max\{\gamma^j(p), 0\}$.
 For each random binary vector $x$ that follows $p$, the following holds due to the union bound.
  $$\Pr(Ax \not \geq b) = \Pr\left(\exists j \in [m] \mid \iprod{A_j, x} < b_j\right) \leq  \sum_{j \in [m]} \Pr\left(\iprod{A_j, x} < b_j\right).$$
  For each $j \in [m]$, let $s > 0$ and $\eta_j = \|A_j\|_2$, we compute
  \begin{align*}
    \Pr(\iprod{A_j, x} < b_j) &= \Pr\left(\mathbb E[\iprod{A_j, x}] -  \iprod{A_j, x} >  \mathbb E[\iprod{A_j, x}] - b_j\right)\\
    &= \Pr\left(\iprod{A_j, p-x} >  \gamma^j(p)\right)\\
    &\leq \exp(- s\gamma^j(p)) \mathbb E[\exp(s\iprod{A_j, p-x})]\\
    &= \exp(- s\gamma^j(p)) \prod_{i \in [n]}\mathbb E[\exp(-sA_{ji}(x_i-p_i))]\\
    &\leq \exp(- s\gamma^j(p)) \prod_{i \in [n]}\exp\left(s^2 A^2_{ji}/8\right)\\
    & = \exp\left(s^2\eta_j^2/ 8 - s\gamma^j(p)\right),
  \end{align*}
  where the first inequality is due to the Chernoff bound $\Pr(X \geq a) \leq \exp(-sa)\mathbb E[\exp(sX)]$ for every $s > 0$, and the second is an application of Hoeffding's lemma. To obtain the tightest bound, we minimize the last expression within the domain $s \geq 0$ to obtain $s^* = 4\gamma^j_+(p)/\eta_j^2$. This provides the bound $\Pr(\iprod{A_j, x} < b_j) \leq \exp(-2[\gamma^j_+(p)]^2/ \eta_j^2)$ for every $j$. Then, we have
 $$\Pr(Ax\not\geq b) \leq m \exp(-2\gamma_+^2(p)/\eta^2) = \exp\left(-2\gamma_+^2(p)/\eta^2 + \log m\right) =\exp\left(-2(\gamma_+^2(p)/\eta^2 - \log m/2)\right).$$
 Note that when $\gamma_+^2(p)/\eta^2 < \log m/2$, this bound becomes trivial. Thus, we can tighten it as $\exp(-2[\gamma_+^2(p)/\eta^2 - \log m/2]_+)$, which implies the claim with $k$ independent samples.
\pfend

\totuni*
\pfstart
\blue{
By partitioning the matrices and vectors according to the index sets \(I\) and \(\bar I\), and substituting \(x_I=s+Sx_{\bar I}\) into \eqref{eq:bip}, we obtain the stated reformulation with the corresponding transformed model parameters, together with the last two constraints enforcing \(x_I\in[0,1]^{|I|}\). Since \(B_{JI}\) is a nonsingular square submatrix of an integral totally unimodular matrix, \(\det(B_{JI})=\pm1\), and hence \(B_{JI}^{-1}\) is integral. Therefore, \(s\) and \(S\) are integral because both $d_J$ and $B_{J \bar I}$ are integral, so binary \(x_{\bar I}\) implies integral \(x_I\). Combined with \(x_I\in[0,1]^{|I|}\), this yields \(x_I\in\{0,1\}^{|I|}\), recovering the only binary requirement omitted from the reformulation.
}
 % By splitting matrices and vectors by the index sets $I$ and replacing $x_I=s + Sx_{\bar I}$ into \eqref{eq:bip}, we obtain the above reformulation with computed new model parameters along with the last two constraints to ensure $x_I \in [0,1]^n$.
 % Compared with the original problem, the only missing constraint is the binary requirement on $x_I$, which is guaranteed since $B_J$ is a TU matrix.
\pfend

\section{Additional Experiment on TSP}
% \section{Mathematical Proofs}
\label{sec:tsp}

Given a set of cities $\mathcal N = \{1,2,\dots,n\}$ with pairwise distances $c$, TSP seeks the shortest Hamiltonian cycle that starts at the depot $1$ and visits every city exactly once. Since the GFORS framework does not support formulations with an exponential number of constraints, the Dantzig–Fulkerson–Johnson (DFJ) formulation \citep{dantzig1954solution} is not applicable. Instead, we adopt the Miller–Tucker–Zemlin (MTZ) formulation \citep{miller1960integer} that contains non-binary integer variables. 
Thus, the TSP problem class is structurally different from the other five classes we tested and requires an additional binary encoding step for reformulation.
We use $\mathcal N_k := \{k,k+1,\dots,n\}$ to denote the subset of cities starting at index $k$:
\[
  \begin{aligned}
    \min_{x \in \{0,1\}^{n^2}} \quad & \iprod{c, x}\\
    \text{s.t.} \quad 
    & \sum_{j \neq i} x_{ij} = 1, && \forall i \in \mathcal N,\\
    & \sum_{i \neq j} x_{ij} = 1, && \forall j \in \mathcal N,\\
    & u_i - u_j + n x_{ij} \leq n-1, && \forall i,j \in \mathcal N_2,\ i \neq j,\\
    & 2 \leq u_i \leq n, && \forall i \in \mathcal N_2.
  \end{aligned}
\]
Since this is not a pure BIP due to the integer variables $u$, we apply the unary (sequential) encoding technique, motivated by its favorable performance on Ising machines \citep{tamura2021performance}. Specifically, we introduce binary variables $y_{ik}$ for $i \in \mathcal N_2$ and $k \in \mathcal N_3$, where $y_{ik}=1$ if and only if $u_i \geq k$. Because $u_i \geq 2$, we always have $y_{i1}=y_{i2}=1$, which implies $u_i = \sum_{k \in \mathcal N_3} y_{ik} + 2$, yielding the following unary-encoding reformulation,
\[
\begin{aligned}
  \min_{\substack{x_{ij} \in \{0,1\} \\ y_{ik} \in \{0,1\}}} \quad 
  & \iprod{c, x} \\
  \text{s.t.} \quad 
  & \sum_{j \ne i} x_{ij} = 1, && \forall i \in \mathcal N,\\
  & \sum_{i \ne j} x_{ij} = 1, && \forall j \in \mathcal N,\\
  & \sum_{k \in \mathcal N_3} y_{ik} - \sum_{k \in \mathcal N_3} y_{jk} + n x_{ij} \leq n-1, && \forall i,j \in \mathcal N_2,\ i \neq j,\\
  & y_{i,k-1} \geq y_{ik}, && \forall i \in \mathcal N_2,\ k \in \mathcal N_4.
\end{aligned}
\]

\begin{table}[!tbp]
  \centering
  \blue{
 \scriptsize
\begin{tabular}{rrrrrrrrrrr}
\toprule
Config & \multirow{2}{*}{NNZ}
& \multicolumn{4}{c}{Avg. Objective Value}
& \multicolumn{4}{c}{Avg. Runtime (sec.)}
& \multirow{2}{*}{Avg. TU Time} \\
\cmidrule(lr){3-6}\cmidrule(lr){7-10}
$n$ & & GRB & cuOpt & TU & No TU
& GRB & cuOpt & TU & No TU & \\
\midrule
$3$ & $1.60e1$ & $14.3$ & $14.3$ & $14.3$ & $14.3$
& $0.00$ & $0.20$ & $1.00$ & $0.87$ & $0.00$ \\
$6$ & $2.65e2$ & $22.5$ & $22.5$ & $22.5$ & $\mathrm{TL}$
& $0.01$ & $0.31$ & $1441.78$ & $1800.00$ & $0.00$ \\
$12$ & $2.76e3$ & $35.7$ & $35.7$ & $\mathrm{TL}$ & $\mathrm{TL}$
& $0.20$ & $0.77$ & $1800.00$ & $1800.01$ & $0.00$ \\
$24$ & $2.48e4$ & $50.9$ & $53.0$ & $\mathrm{TL}$ & $\mathrm{TL}$
& $2.28$ & $746.15$ & $1800.01$ & $1800.01$ & $0.00$ \\
$48$ & $2.10e5$ & $82.7$ & $149.5$ & $\mathrm{TL}$ & $\mathrm{TL}$
& $56.68$ & $1800.49$ & $1800.01$ & $1800.01$ & $0.01$ \\
$96$ & $1.72e6$ & $173.6$ & $\mathrm{TL}$ & $\mathrm{TL}$ & $\mathrm{TL}$
& $1386.09$ & $1800.49$ & $1800.01$ & $1800.01$ & $0.05$ \\
$192$ & $1.40e7$ & $945.8$ & $\mathrm{TL}$ & $\mathrm{TL}$ & $\mathrm{TL}$
& $1272.88$ & $1800.61$ & $1800.01$ & $1800.01$ & $0.55$ \\
$384$ & $1.13e8$ & $\mathrm{TL}$ & $\mathrm{TL}$ & $\mathrm{TL}$ & $\mathrm{TL}$
& $1806.41$ & $2095.93$ & $1800.04$ & $1800.02$ & $5.06$ \\
$768$ & $9.03e8$ & $\mathrm{OOM}$ & $\mathrm{OOM}$ & $\mathrm{OOM}$ & $\mathrm{OOM}$
& $1055.41$ & $8.08$ & -- & -- & -- \\
\bottomrule
\end{tabular}

 }
 \caption{TSP: objective values and runtimes. \emph{TU} and \emph{No TU} denote the implementations of GFORS \blue{with and without} TU reformulation, respectively; \blue{$\mathrm{TL}$ and $\mathrm{OOM}$ indicate no feasible solution before the time or memory limit, respectively.}}
  \label{tab:tsp}
\end{table}

We consider sizes $n$ ranging from $3$ to $768$, \blue{with the symmetric distance matrix $C$ generated using entries sampled uniformly from $\{1,\ldots,9\}$ and with zero diagonal entries.}
Results are reported in Table~\ref{tab:tsp}. 
In this case, all equality constraints can be eliminated via TU reformulation. The columns \emph{TU} and \emph{No TU} correspond to implementations with and without this step. However, both GFORS variants perform unsatisfactorily: \blue{the TU variant solves only the two smallest instances, while the No-TU variant solves only the smallest one. 
  In comparison, GRB remains effective through $n=192$, while cuOpt matches GRB through $n=12$ but deteriorates thereafter, reaching the time limit from $n=96$ onward and running out of memory at $n=768$.} This suggests that the loose relaxation of the MTZ formulation and the extremely sparse feasible set relative to the domain $[0,1]^n$ are potentially the main bottlenecks for GFORS in this setting.

% We consider sizes $n$ ranging from $3$ to $768$, \blue{with the symmetric distance matrix $C$ generated using entries sampled uniformly from $\{1,\ldots,9\}$ and with zero diagonal entries.}
% Results are reported in Table~\ref{tab:tsp}. 
% In this case, all equality constraints can be eliminated via TU reformulation. The columns \emph{TU} and \emph{No TU} correspond to implementations with and without this step. However, both GFORS variants perform unsatisfactorily: \blue{only the smallest instance is solved, while GRB remains effective through $n=192$ and cuOpt through $n=48$ on the same formulation.} This suggests that the loose relaxation of the MTZ formulation and the extremely sparse feasible set relative to the domain $[0,1]^n$ are potentially the main bottlenecks for GFORS in this setting.

% In this case, all equality constraints can be eliminated via TU reformulation. The columns \emph{TU} and \emph{No TU} correspond to implementations with and without this step. However, both variants perform unsatisfactorily: only the smallest instances are solved, while Gurobi can handle sizes up to $n=96$ for the same formulation. This suggests that the loose relaxation of the MTZ formulation and the extremely sparse feasible set relative to the domain $[0,1]^n$ are potentially the main bottlenecks for GFORS in this setting.

\end{document}